\documentclass[leqno]{article}
\usepackage[T1]{fontenc}
\usepackage{upgreek}
\usepackage{amsmath,amsfonts,amssymb,amsthm,mathtools,latexsym}
\usepackage{eucal}
\usepackage{indentfirst}
\usepackage{tikz}
\usetikzlibrary{calc}
\usetikzlibrary{arrows,decorations.pathmorphing,backgrounds,positioning,fit,petri,arrows.meta} 
\usepackage[scr]{rsfso}
\usepackage{graphics}
\usepackage{bbm} 
\usepackage{anyfontsize}
\usepackage{stmaryrd}
\usepackage{mbboard/texinputs/mbboard} 
\newcommand{\sigmapqtimes}{\underset{\Sigma_p\times \Sigma_q}{\times}}
\newcommand{\sigmaqpt}{\underset{\Sigma_q\times \Sigma_p}{\times}}
\newcommand{\pre}{\text{\normalfont{pre-}}}
\newcommand{\myhom}{\text{\normalfont{Hom}}}
				
\newcommand{\xtworightarrow}[2][]{%
				  \xrightarrow[#1]{#2}\mathrel{\mkern-14mu}\rightarrow
	}

\newcommand\rightTwoArrows{%
        \mathrel{\vcenter{\mathsurround0pt
                \ialign{##\crcr
												\noalign{\nointerlineskip}$\xrightarrow{\qquad}$\crcr
												\noalign{\nointerlineskip}$\xrightarrow{\qquad}$\crcr
                }%
        }}%
}
\newcommand\mapstoto{%
        \mathrel{\vcenter{\mathsurround0pt
                \ialign{##\crcr
												\noalign{\nointerlineskip}$\xmapsto{\qquad}$\crcr
												\noalign{\nointerlineskip}$\xmapsto{\qquad}$\crcr
                }%
        }}%
}
\newcommand{\plusideal}{\mathfrak{I}} 
\newcommand{\plusprimeideal}{\mathfrak{p}} 
\newcommand{\step}{\item[(\addtocounter{equation}{1}\theequation) ]}
\newcommand{\op}{\text{op}}
\newcommand{\Bnabla}{\bbnabla}
\newcommand{\BDelta}{\bbDelta}
\newcommand{\asetnab}{\text{$A$\normalfont-Set}^\Bnabla}
\newcommand{\aset}{\text{$A$-Set}}

\newcommand{\bset}{\text{$B$-Set}}
\makeatletter
\newcommand*{\rom}[1]{\expandafter\@slowromancap\romannumeral #1@}
\makeatother
\newcommand{\cring}{\text{CRing}}
\newcommand{\characteristic}{{\mathbf{\mathbbm{1}}}}
\newcommand{\unit}{{\mathbf{1}}}
\newcommand{\bigo}{\mathcal{O} }
\newcommand{\gr}{\text{GR}}

\newcommand{\cgr}{\text{CGR}}
\newcommand{\tgr}{\text{TGR}}
\newcommand{\set}{\text{Set}}
\newcommand{\Q}{\mathbb{Q}}
\newcommand{\R}{\mathbb{R}}
\newcommand{\C}{\mathbb{C}}
\newcommand{\N}{\mathbb{N}}

\newcommand{\fin}{\mbox{Fin}}
\newcommand{\cok}{\mbox{Cok}}
\newcommand{\kernel}{\mbox{Ker}}
\newcommand{\cokernel}{\mbox{CokKer}}
\makeatletter
\DeclareRobustCommand{\properideal}{\mathrel{\text{$\m@th\proper@ideal$}}}
\newcommand{\proper@ideal}{%
  \ooalign{$\lneq$\cr\raise.22ex\hbox{$\lhd$}\cr}%
}
\makeatother
\newcommand{\lideal}{\lhd}

\newcommand{\Z}{\mathbb{Z}}
\newcommand{\ch}{\text{Ch}}
\newcommand{\D}{\mathcal{D}}

\newcommand{\E}{\mathcal{E}}
\newcommand{\crg}{\text{CRig}}
\newcommand{\mymod}{\text{\normalfont-mod}}

\newtheorem{definition}[equation]{Definition} 
\newtheorem{definition*}{Definition} 
\newtheorem{remark}[equation]{Remark}
\newtheorem*{remark*}{Remark}

\numberwithin{equation}{section}
\title{Algebra over generalized rings}
\author{Shai Haran \\ {haran@technion.ac.il} }
\date{}
\begin{document}
\maketitle
\begin{abstract}
				For a commutative ring $A$, we have the category of (bounded-below) chain
				complexes of $A$-modules $Ch_{+}(A\mymod)$, a closed symmetric monoidal
				category with a compatible stable Quillen model structure. The associated homotopy
				category is the derived category $\mathbbm{D}(A\mymod)$, where one inverts all the
				quasi-isomorphisms, and it has the good description as 
				\begin{equation}
				\mathbbm{D}(A\mymod)=\ch_{+}(\text{proj-}A\mymod)\,/\simeq
												\label{eq0.1}
				\end{equation}
				the chain complexes made up of projective $A$-module in each dimension, and chain
				maps taken up to chain homotopy. We give here the analogous theory for a
				(commutative) generalized ring in the sense of \cite{MR3605614}. We use here a
				different concept of ``A-module'' than the one used in \cite{MR3605614} (which was
								useful for the introduction of derivatives and the cotangent complex in the
				arithmetical settings). We refer to the new concept as ``$\aset$''. For an ordinary commutative ring $A$, an $A$-set is just
				an $A$-module in the usual meaning, and our construction will be equivalent to
				$\mathbbm{D}(A\mymod)$. For the initial object of the category of generalized rings
				$\mathbb{F}$ ``the field with one element'', we obtain the category of symmetric
				spectra, and the associated stable homotopy category with its smash product (an
				$\mathbb{F}$-set is just a pointed set, i.e. a set $X$ with a distinguish element
$O_X\in X$). Thus the analogous theories of stable homotopy and of chain
				complexes of modules over a commutative ring appear as two sides of the same coin,
				and moreover, they appear in a context where they interact (via the forgetful
				functor and its left adjoint - the base change functor). For the ``real
				integers'' $A=\Z_{\R}$, the $\Z_{\R}$-sets include
				the symmetric convex subsets of $\R$-vector spaces. We also give the global theory
				of the derived category of $\bigo_X$-sets, for a generalized scheme $X$, in a way
				that is based on the local projective model structure rather than follow the path
				of Grothendieck \cite{grothendieck1961elements} Jardine \cite{MR906403}  Voevodsky \cite{MR1648048}, of using
				injective resolution (but inevitably uses $\infty$-categories). 
\end{abstract}
\section{Generalized Rings cf. \cite{MR3605614}}
Let  $\fin_{0}$ denote the category of finite pointed sets. For a map
$f\in\fin_{0}(X,Y)$, we have its kernel and cokernel 
\begin{equation}
				\begin{array}[H]{cl}
								(\romannumeral 1) & \kernel(f): f^{-1}(O_Y) \rightarrow X \\
								(\romannumeral 2) & \cok(f): Y \twoheadrightarrow Y/f(X)
				\end{array}
				\label{eq:1.1}
\end{equation}
where $Y/f(X)$ is obtained from $Y$ by contracting $f(X)$ to a point. There is a
canonical map 
\begin{equation}
				\tilde{f}: X/f^{-1}(O_Y)= \cokernel(f) \longrightarrow \kernel\cok (f) =f(X)
				\label{eq:1.2}
\end{equation}
which is not always an isomorphism, and we let $\mathcal{F}_0$ denote the subcategory of
$\fin_{0}$ with the same objects, but with the maps $f$ for which $\tilde{f}$ is an
isomorphism, $X/f^{-1}(O_Y) \cong f(X)$, that is $f$ is bijective away from $O_Y$:
\begin{equation}
				f(x_1)=f(x_2) \not= O_Y \Longrightarrow x_1=x_2. 
				\label{eq:1.3}
\end{equation}
We view $\mathcal{F}_0$ as the category of ``finite dimensional vector spaces over the
field with one element''. We will omit the pointed element $O_X\in X$ for $X\in \fin_0$
and we have isomorphic category $\fin_\bullet$, where the objects are finite sets and the maps
are partially-defined functions 
\begin{equation}
				\begin{array}[H]{rl}
								\fin_0 &\xleftrightarrow{\;\; \sim \;\;} \fin_\bullet \\
				X & \longmapsto X\setminus \{O_X\} \\
				n_{+} := n \coprod\left\{ O_n \right\}  &\longmapsfrom n
				\end{array}
				\label{eq:1.4}
\end{equation}
We denote by $\mathcal{F}\subseteq \fin_\bullet$ the subcategory corresponding to
$\mathcal{F}_0\subseteq \fin_0$. Thus $\mathcal{F}$ has objects the finite sets, and the
maps are the partially defined bijections 
\begin{equation}
				\mathcal{F}(m,n) = \left\{ f \; : \; D(f)\xrightarrow{\sim} I(f) , \quad
				D(f)\subseteq m, \quad I(f)\subseteq n \right\}
				\label{eq:1.5}
\end{equation}
Note that the category $\mathcal{F}$ is self-dual, $\mathcal{F}\simeq \mathcal{F}^{\op}$, and
we have an involution \linebreak $f\mapsto f^t$,
$\mathcal{F}(m,n)\xrightarrow{\sim}\mathcal{F}(n,m)$, taking $f=\left(
f:D(f)\xrightarrow{\sim} I(f) \right)$ to  \linebreak $f^{t}=\left( f^{-1}:I(f)\xrightarrow{\sim}
D(f) \right)$,  satisfying 
\begin{equation}
				f^{tt}=f , \quad \left( f_2\circ f_1 \right)^{t}=f_1^t\circ f_2^t, \quad
				(\text{Id}_n)^t = \text{Id}_n. 
				\label{eq:1.6}
\end{equation}
The initial and final object of $\mathcal{F}$ is thus the empty set $\emptyset$.
\begin{itemize}
				\step 
								{\bf{Definition}}\cite{MR3605614}:
 \begin{textit} 
				 A generalized ring is a fuctor $A:\mathcal{F}\to \text{Set}_0$, pointed
				 \linebreak
				$A_{\emptyset}=\left\{ 0 \right\}$, and with operations of ``multiplication'' and
				``contraction'' as follows. For $f\in \fin_{\bullet}(m,n)$ put
				$A_f=\prod\limits_{j\in n}A_{f^{-1}(j)}$. For $b=(b_j)\in A_f$, $b_j\in
				A_{f^{-1}(j)}$, we usually omit $f$ and write ``$m\xrightarrow{b}n$''.
				\begin{itemize}
								\item[(\romannumeral 1)] {\bf{multiplication}}: $
												\begin{array}[t]{ll}
												A_n\times A_f\to A_m \\
												a, (b_j) \mapsto a\lideal b
												\end{array} $
								\item[(\romannumeral 2)] {\bf{contraction}}: 
												$\begin{array}[t]{ll}
												A_m\times A_f\to A_n \\
												c, (b_j) \mapsto c\sslash b
												\end{array} $
				\end{itemize}
				\label{NFG}
\end{textit}
For a map $g\in \fin_{\bullet}(n,q)$, we can extend these maps fiberwise  to get       
				\step
				\begin{itemize}
								\item[(\romannumeral 1) ] {\bf{multiplication}}: $
												\begin{array}[t]{ll}
																A_g\times A_f \to A_{g\circ f} \\
																a,b \mapsto a\lideal b, \quad (a\lideal b)_i = a_i \lideal
																(b_j)_{f(j)=i} \quad \text{for $i\in q$} 
												\end{array} 
												$
								\item[(\romannumeral 2)] {\bf{Contraction}}: 	 $
												\begin{array}[t]{ll}
																A_{g\circ f}\times A_f\to A_g \\
																c,b \mapsto c\sslash b, \quad (c\sslash b)_i =
																c_i\sslash(b_j)_{f(j)=i}, \quad i\in q
												\end{array}
												$ 
				\end{itemize}
				We require these operations to satisfy the following axioms. 
				\step \label{1.9} {\underline{\bf{Associativity}}}: $(a\lideal b)\lideal e = a\lideal(b\lideal e)$
				\step \label{1.10} {\underline{\bf{Unit}}:} We have $\unit\in A_{[\unit]}$, (with [1] denoting the
				set with one element):
				\[
								\begin{array}[H]{ll}
												\unit \lideal a = a \\
												a \lideal (\unit)_{j\in n} = a  = a\sslash (1)_{j\in n}, \qquad
												a\in A_n.
								\end{array}
				\]
\end{itemize}
				Thus, forgetting the contraction, $A$ is just an operad. We require the
				following identities to hold (in $A_{k\to q}$, where we can take $q=\left\{
								\ast
\right\}=[\unit])$ 
				\begin{itemize}
				\step {\underline{\bf{Left-Adjunction:}}} $(d\sslash c)\sslash a =
								d\sslash (a\lideal c)$ \\
								\begin{center}
												\begin{tikzpicture}
																\node (n) at (4,2) {$n$};
																\node (k) at (5.5,0) {$k$};
																\node (m) at (4,0.8) {$m$};
																\node (q) at (3,-1) {$q$};
																\draw [->] (n) edge[post,bend right,left] node{$d$} (q);
																\draw [->] (n) edge[right] node{$c$} (m);
																\draw [->] (n) edge [post,bend left,right]  node{$a\lideal c$} (k);
																\draw [->] (m) edge[right,above] node{$a$} (k);
																\draw [->] (m)-- (q);
																\Huge 
																				\draw[double,thick,-{Implies[]}] (k) --++(q);
												\end{tikzpicture}
								\end{center}
								\step {\underline{\bf{Right-Adjunction:}}} $d\sslash (a\sslash c)
								= (d\lideal c)\sslash a$ \\
								\begin{center}
																\begin{tikzpicture}
																				\node (n) at (4,2) {$n$};
																				\node (k) at (5.5,0) {$k$};
																				\node (m) at (4,0.8) {$m$};
																				\node (q) at (3,-1) {$q$};
																				\draw [->] (m) edge[post,left] node{$d$} (q);
																				\draw [->] (n) edge[right] node{$c$} (m);
																				\draw [->] (n) edge [post,bend left,right]  node{$a$} (k);
																				\draw [->] (m) edge[right,above]
																				node{$\scriptstyle a\sslash c$} (k);
																				\Huge 
																								\draw[double,thick,-{Implies[]}] (k) --++(q);
																\end{tikzpicture}
								\end{center}
								\step {\underline{\bf{Left-Linear:}}} $d\lideal(a\sslash c) = (d\lideal a)\sslash c $
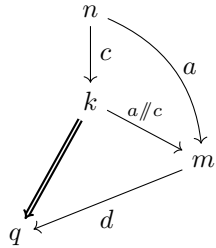
								\begin{center}
																\begin{tikzpicture}
																				\node (n) at (4,2) {$n$};
																				\node (m) at (5.5,0) {$m$};
																				\node (k) at (4,0.8) {$k$};
																				\node (q) at (3,-1) {$q$};
																				\draw [->] (m) edge[post,below] node{$d$} (q);
																				\draw [->] (n) edge[right] node{$c$} (k);
																				\draw [->] (n) edge [post,bend left,right]  node{$a$} (m);
																				\draw [->] (k) edge[right,above]
																								node{$\scriptstyle a\sslash c$} (m);
																				\Huge 
																								\draw[double,thick,-{Implies[]}] (k) --++(q);
																\end{tikzpicture}
								\end{center}

				\step {\underline{\bf{Right-Linear:}}} $(d\sslash c)\lideal a = (d\lideal
								\tilde{a})\sslash \tilde{c} $
		
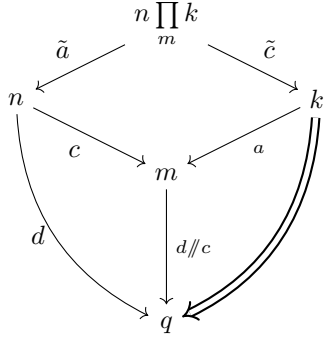
								\begin{center}
																\begin{tikzpicture}
																				\node (nmk) at (4,2) {$n\prod\limits_{m}k$};
																				\node (n) at (2,1) {$n$};
																				\node (k) at ((6,1) {$k$};
																				\node (m) at (4,0) {$m$};
																				\node (q) at (4,-2) {$q$};
																				\draw [->] (nmk) edge[post,above left] node{$\tilde{a}$} (n);
																				\draw [->] (m) edge[post,right] node{$\scriptstyle d\sslash c$} (q);
																				\draw [->] (n) edge[below left] node{$c$} (m);
																				\draw [->] (nmk) edge [post,above right]
																								node{$\tilde{c}$} (k);
																				\draw [->] (n) edge [post,bend right,left]  node{$d$} (q);
																				\draw [->] (k) edge[right,below right]
																								node{$\scriptstyle a$} (m);
																				\draw[double,thick,double distance=2pt]
																								(k) edge[double,thick,-{Implies[]},bend left] (q);
																\end{tikzpicture}
								\end{center}
								\end{itemize}
								Here $\tilde{a}$ (resp. $\tilde{c}$) is obtained from $a$ (resp.
								$c$) by identifying the fibers; i.e. denoting by $f:n\to m$ we have
								\begin{equation}
												\tilde{a}=(\tilde{a}_j) \text{ with }\tilde{a}_j=a_{f(j)}
												\label{eq:1.15}
								\end{equation}
								It follows from the axioms that for a generalized ring $A$, the set
								$A_{[1]}$ is an associative (1.9), unital (1.10), but also commutative
								(1.14), monoid with respect to multiplication. It has the involution
								\begin{equation}
												a^t := 1\sslash a, \qquad a\in A_{[1]}. 
												\label{eq:1.16}
								\end{equation}
												We denote by $A^{+}_{[1]}$ the sub-monoid of symmetric element. 
								\begin{equation}
												A^{+}_{[1]}=\left\{ a\in A_{[1]}, a^t=a \right\}
												\label{eq:1.17}
								\end{equation}
								Note that the axioms (1.11-14) imply that any formula made up of multiplication
								$\lideal$, and contraction $\sslash \;$,  is equivalent to a formula with
												{\underline{one}} contraction of the form 
												\begin{equation}
																(a_1\lideal a_2\lideal \dots \lideal a_m)\sslash (b_1\lideal b_2\lideal \dots \lideal b_n)
																\label{eq:1.18}
												\end{equation}
												A homomorphism of generalized rings $\varphi: B\to A$ is a natural
												transformation $\varphi_n: B_n\to A_n$ preserving multiplication,
												contraction, and the unit $1$ (the unit $0$ is always preserved since for
												$n=\emptyset$, $B_{\phi}=\left\{ 0 \right\}= A_{\phi}$). We thus have a category of generalized rings $\gr$
												(these are what we called ``commutative-generalized-ring`` in
												\cite{MR3605614}, because the axiom (1.14) imply the commutativity of
								$A_{[1]}$ ). \vspace{.1cm}\\
												The category $\gr$ is complete and co-complete. Limits and filtered co-limits are
												formed in $\text{Set}_0$. More general co-limits are more complicated, but we do
												have push-outs
								\begin{equation}
								\begin{array}[h]{l}
																\begin{tikzpicture}
																				\node (A1) at (4,2) {$A_1$};
																				\node (B) at (2,2) {$B$};
																				\node (A0) at (2,0) {$A_0$};
																				\node (A01) at (4.4,0) {$A_0\underset{B}{\otimes}A_1$};
																				\node (A11) at (4,0.2) {};
																				\draw[->] (B)--(A1);
																				\draw[->] (B)--(A0);
																				\draw[->] (A0)--(A01);
																				\draw[->] (A1)--(A11);
																				\draw[-] (3.9,1)--(3,1);
																				\draw[-] (3,1)--(3,0.1);
																\end{tikzpicture}
								\end{array}
								\label{eq:1.19}
								\end{equation}
								and in particular (categorical) sums
								$A_0\underset{\mathbb{F}}{\otimes}A_1$. \\
								We let $\tgr\subseteq \gr$ denote
								the full subcategory with objects the generalized rings $A$ satisfying the
								''total-commutativity``, using the notations of (1.15). 
								\begin{equation}
								c\lideal \tilde{a}= a\lideal \tilde{c} \quad \text{ in }  \quad  A_{n\underset{m}{\prod} k\to m}
												\label{eq:1.20}
								\end{equation}
								(we can take $m=\left\{ \ast \right\}=[1]$, so $a\in A_k$, $c\in A_n$,
								$c\lideal \tilde{a}= a\lideal \tilde{c}\in A_{n\prod k}$). \\
								The inclusion
								$\tgr \hookrightarrow \gr$ has a left-adjoint, $A\mapsto A^{T}$, where $A^{T}$ is
								the maximal totally commutative quotient of $A$. \vspace{.1cm}\\ 
								Here are some examples of totally-commutative generalized rings. 
				\begin{itemize}
				\step $\underline{\mathbb{F}}$. \\
								The initial object of $\gr$, ''the field with one element``, is denoted by
								$\mathbb{F}$, $\mathbb{F}: \mathcal{F}\xrightarrow{\sim}
								\mathcal{F}_0\subseteq \text{Set}_0$, $\mathbb{F}_n := n \coprod \left\{
								0_n \right\}=\left\{ \delta_{i} \right\}_{i\in n}\coprod \left\{ 0_n
								\right\} $
				\step {\underline{$\mathbb{F}\left\{ M \right\}$}}. \\
								For a commutative monoid $M$, we have $\mathbb{F}\left\{ M \right\}\in
								\tgr$ 
								\[ \mathbb{F}\left\{ M \right\}_{n} := (n\times M) \coprod \left\{ 0_n
								\right\} \]
								This give a full and faithful functor $M\mapsto \mathbb{F}\left\{ M
								\right\}: \text{CMon}\hookrightarrow \tgr $
				\step A {\underline{Commutative-rig}} (or a ''semi-ring``) is a set
								$A$ with two associative and commutative operations of addition $+$, with
								unit $0$, and multiplication $\cdot$, with unit $1$, where $x\cdot 0 = 0 $
								for all $x\in A$, and we have distributivity $x\cdot (y_1+y_2)=(x\cdot
								y_1)+(x\cdot y_2)$. These form a category $\crg$, where the arrows are
								the set maps preserving the operations and the units $0,1$. For
								$A\in\crg$, we have $A\in\tgr$ (denoted by the same letter!), with
								$A_n:=A^n$, and with the operations of multiplication and contratction
								defined as follows:
								\[
												\text{Note that}\quad A_f \equiv \prod\limits_{j\in n} A^{f^{-1}(j)} =
								A^{\coprod\limits_{j\in n}f^{-1}(j)} = A^{D(f)}, \]
				\begin{itemize}
								\item[(\romannumeral 1) ]
				 $ 
				 	\begin{array}[t]{ll}
				 \lideal : A_{n}\times A_{f}\to A_m \\\\
								(a_j)_{j\in n}, \quad (b_i)_{i\in D(f)}\mapsto (a\lideal b)_i &:=
								a_{f(i)}\cdot b_i, \quad i \in D(f), \\
								&:= 0, \quad i\in m\setminus D(f) 
				\end{array}
								$
				\item[(\romannumeral 2) ] 
				 $ 
				 	\begin{array}[t]{ll}
				 \sslash : A_{m}\times A_{f}\to A_n \\\\
				 (c_i)_{i\in m}, \; (b_i)_{i\in D(f)}\mapsto (c\sslash b)_j &:=
				 \sum\limits_{i\in f^{-1}(j)} c_i\cdot b_i, \quad j\in I(f) \\
				 & := 0, \quad j\in n\setminus I(f)
 \end{array}
				 $
				\end{itemize}
				This gives a full and faithfull embedding $\crg\hookrightarrow \tgr$. \\ Examples of
				commutative-rigs, include commutative-rings, but there are many more examples where addition
				is not invertible, such as the ''tropical`` rigs 
								 \[ \left\{ 0,1 \right\}\hookrightarrow
								 [0,1]\hookrightarrow [0,\infty)\equiv \R^{+}_{\max} \]
				with the usual multiplication, and with $y_1+y_2:=\max \left\{ y_1,y_2 \right\}$.
\step {\underline{$\Z_\R$ and $\Z_\C$}} \\
				We have the (maximal compact sub-topological-generalized ring), the ''real
				integers`` $\Z_\R\subseteq \R$ (resp. the ''complex integers`` $\Z_\C\subseteq
				\C$), given by 
				\[ (\Z_{\R})_{n}:=\left\{ (a_j)\in\R^n,\;\; \sum\limits_{j\in n}|a_j|^{2}\le 1
				\right\} \]
				resp. 
				\[ (\Z_\C)_n := \left\{ (a_j)\in \C^n,\;\; \sum\limits_{j\in n}|a_j|^2\le 1
				\right\} \] 
				the unit $\ell_2-\text{ball}$ in $\R^n$, resp $\C^n$.
				\end{itemize}
				\begin{remark}
								 We prefer to work with $\gr$, rather than $\tgr$,
								since the categorical sum in $\gr$ of $\Z$ with itself, $\Z\otimes_\mathbb{F}\Z$, is a
								very interesting object, while the categorical sum in $\tgr$ reduces to the
								ordinary integer: $(\Z\otimes_\mathbb{F}\Z)^{T}\xrightarrow{\sim}\Z$ (and so the
								''Arithmetical surface`` reduces to its diagonal, as it does with ordinary
								rings). 
								\label{rem:1.25}
				\end{remark}
				\begin{definition}
We define the ''{\underline{commutative}}``
				generalized rings $\cgr$ to be the intermediate full sucategory
				\[ \tgr \subseteq \cgr \subseteq \gr \]
				where $A\in \gr$ is \underline{commutative} if for $a,b\in A_n$, $c,d\in A_m$ 
				$(x_{i,j})\in\left( A_{[1]} \right)^{m\prod n}$, we have (in $A_{[1]}$): 
\begin{equation}
				  (a\lideal \tilde{d}\lideal (x_{i,j}))\sslash (b\lideal\tilde{c})=(d\lideal
												\tilde{a}\lideal (x_{i,j}))\sslash (c\lideal \tilde{b}) 
												\label{eq:1.27}
\end{equation}
\end{definition}
\section{$A$-sets}
\begin{definition}
				\itshape	For $A\in \cgr$, an $A$-set is a pointed set $M\in\set_0$ together with
				an $A$-action: For each $n\in \mathcal{F}$ we have maps
				\[
								\begin{array}[H]{rl}
												A_n\times M^n\times A_n &\xrightarrow{\quad} M \\
												(b,m_j,d)&\xmapsto{\quad} \langle b,m,d\rangle_n
								\end{array}
				\]
				These maps are required to satisfy the following axioms:
				\begin{description}
								\item[{\underline{Associativity}}: ] For a map $f\in\set_\bullet(m,n)$, for
												$b,d\in A_n$, $b^{\prime},d^{\prime}\in A_f$, 
												\[ \langle b\lideal b^{\prime},m_i,d\lideal d^{\prime}\rangle_m=
																\langle b, \langle b^{\prime}_{j},
												m_i,d^{\prime}_{j}\rangle_{f^{-1}(j)},d\rangle_{n}  \] 
								\item[ {\underline{Unit}}: ] $\langle 1,m,1\rangle_{[1]}=m$ \\
												(i.e. we have an action of the operad $A\times A$). 
								\item[ {\underline{naturality}}: ]  For $f\in \set_{\bullet}(m,n)$, $d\in
												A_n$, $b\in A_m$, $a\in A_f$, $(m_j)\in M^{n}$, 
								\[ 
												\begin{array}[h]{ll}
																\langle b,m_{f(i)},d\lideal a\rangle_m = \langle b\sslash a,
																m_j,d\rangle_n \\
																\langle d\lideal a, m_{f(i)}, b \rangle_m = \langle d, m_j,
																b\sslash a\rangle_n 
								\end{array} \]
				\item[{\underline{Commutativity}}: ]For $b,d\in A_n$, $b^{\prime},d^{\prime}\in A_m$, $m_{j,i}\in
								M^{n\times m}$, 
								\[ \langle d\lideal \tilde{d}^{\prime},m_{j,i}, b\lideal
												\tilde{b}^{\prime}\rangle_{n\times m}=
												\langle d^{\prime}\lideal \tilde{d},m_{j,i},b^{\prime}\lideal
								\tilde{b}\rangle_{m\times n}  \]  
				\end{description}
				\end{definition}
								A map $\varphi: M\to N$ of $\aset$ is a set map $\varphi \in
								\set_0(M,N)$, preserving the $A$-action, 
								\[ \varphi \left( \langle b,m_j,d\rangle_n \right) = \langle
								b,\varphi(m_j),d\rangle_n \]
								Thus we have the category $\aset$.
\subsection*{Examples of $\aset$}
\begin{itemize}
				\step Given a homomorphism $\varphi\in\cgr(A,B)$, $B_{[1]}$ is an
								$\aset$: 
								\[ \langle b,x_i,d\rangle_n := \left( \varphi(b)\lideal (x_i) \right)\sslash
								\varphi(d), \qquad b,d\in A_n, x_i\in \left( B_{[1]} \right)^n. \]
								\step A sub-$A$-set of $A_{[1]}$ is just an ideal, cf.\cite{MR3605614}.
								\step For $A=\mathbb{F}$ we have $\mathbb{F}\text{-set}\equiv \set_0$, and for
								$M\in \set_0$ there is a unique $\mathbb{F}$-action:
								\[ \langle b,m_j,d\rangle_n = \left\{
												\begin{array}[h]{ll}
																m_i & \text{if $b=d=\delta_i$} \\
																0 & \text{otherwise}
				\end{array}\right. \]
\step For a commutative ring $A\in \cring\subseteq \cgr$, we have
				\[\aset\equiv A\mymod\]
				the category of $A\text{\normalfont-modules}$.
\step For the real (resp. complex) integers $A=\Z_\R$ (resp. $\Z_\C$), the
''finite-dimensional-torsion-free`` $A\text{\normalfont -Sets}$ are the convex subsets
				$M\subseteq \R^n$ (resp. $M\subseteq \C^n$) which are symmetric $u\cdot M\subseteq
				M$ for $|u|\le 1$.
\step {\bf{\underline{Notation:}}}: For $M\in\aset$, and for $m\in M$,
				$a\in A_{[1]}$ we write: 
\[ 
				\begin{array}[h]{ll}
								a\cdot m := \langle a,m,1\rangle_{[1]} \\
								m\cdot a := \langle 1,m,a\rangle_{[1]} = a^t\cdot m, \quad a^t=1\sslash a 
				
\end{array}\]
We have 
\[ a_1\cdot (a_2\cdot m) = (a_1\lideal a_2)\cdot m, \quad 1\cdot m = m , \qquad (a_1\cdot
m)\cdot a_2 = a_1\cdot (m\cdot a_2) \]
\end{itemize}
The category $\aset$ is complete and co-complete. Inverse limits, and filtered
co-limits are formed in $\set_0$, while co-limits and sums are more complicated. Given a
set $V$ the free $\aset$ on $V$, $A^V$, is given by 
\begin{equation}
				A^V = \left( \coprod\limits_{n\in\mathbb{F}} A_n\times V^n\times A_n \right)
				\Big\slash \sim
\end{equation} 
where $\sim$ is the equivalence relation generated by naturality and commutativity.
The element of $A^V$ can be written, non-uniquely, as 
\[ \langle b,v,d\rangle_n, \quad b,d\in A_n, \quad v\in V^n. \]
For $M\in \aset$, we have 
\begin{equation}
				\begin{array}{l}
								\set(V,M)\equiv \aset(A^V,M) \\\\
								\varphi \mapsto \tilde{\varphi}\left( \langle b,v,d\rangle_n \right) = \langle
								b,\varphi(v),d\rangle_n. 
				\end{array}
\end{equation}
When $V=\left\{ v_0 \right\}$ is a singleton, $A^{\left\{ v_0 \right\}}\equiv
A_{[1]}\cdot v_0$, and the free $\aset$ on one generator is just $A_{[1]}$. \\
Given a homomorphism $\varphi\in \cgr(B,A)$, we have an adjunction (we use ''geometric``
notation): 
								\begin{gather}
												\begin{array}[H]{l}
												\begin{tikzpicture}
																\node (a) at (4,4) {$\text{$A$-set}$};
																\node (b) at (4,2) {$\text{$B$-set}$};
																\draw [->] (a) to [out=330,in=30,right] node{$\varphi_{*}$} (b);
																\draw [<-] (a) to [out=210,in=150,left] node{$\varphi^{*}$} (b);
												\end{tikzpicture} 
												\label{eq:2.10}
												\end{array}
								\end{gather}
								where $\varphi_* N\equiv N$ with the $B$-action $\langle
								b,m_j,b^{\prime}\rangle_n := \langle
								\varphi(b),m_j,\varphi(b^{\prime})\rangle_n$. The left adjoint is
								$\varphi^*M=\left( \coprod\limits_{n} A_n\times M^n\times A_n
								\right)\big/\sim$, where $\sim$ is the equivalence relation generated by
								naturality, commutativity, and $B$-linearity: 
								\[ \langle a, \langle
												b_{j},m_i,b^{\prime}_{j}\rangle_{f^{-1}(j)},a^{\prime}\rangle_{n} =
								\langle a\lideal \varphi(b),m_i,a^{\prime}\lideal \varphi(b^{\prime})\rangle_m \] 
								for $f\in \set(m,n)$, $a,a^{\prime}\in A_{n},b,b^{\prime}\in B_f$. \\
								In particular, for $B=\mathbb{F}$ and $\phi\in
								\cgr(\mathbb{F},A)$ the unique homomorphism, $\phi_*$ is just
								the functor forgetting the $A$-action, and for $V\in \set_0$
								$\phi^{*} V= A^{V\setminus \langle 0 \rangle}$ is the free
								$A$-set on $V\setminus \left\{ 0 \right\}$. \vspace{.2cm}\\
								For $M,N,K\in \text{$A$-set}$, let the ''bilinear maps`` be defined by 
								\begin{equation}
												\text{Bil}_A(M,N;K) =   \left\{ \begin{array}[h]{ll} 
																				\varphi: M\wedge N\to K, &   \varphi\left( \langle
																				a,m_j,a^{\prime}\rangle ,n \right) = \langle a, 
																				\varphi(m_j\wedge n),a^{\prime}\rangle, \\ &  \varphi\left(
																								m,\langle a,n_j,a^{\prime}\rangle
												\right)= \langle a,\varphi(m\wedge n_j),a^{\prime}\rangle \end{array}\right\} 
												\label{eq:2.11}
								\end{equation}
								It is a functor in $K$, and as such it is representable 
								\begin{equation}
												\text{Bil}_A(M,N;K)\equiv \text{$A$-set}(M\otimes_A N,K)
												\label{eq:2.12}
								\end{equation}
								Where $M\otimes_A N$ is the free $\aset$ on $M\wedge N$ modulo the
								equivalance relation generated by the $A$-bilinear relations, and where
								$\otimes:M\wedge N\to M\otimes_A N$ is the universal $A$-bilinear map. The
								elements of $M\otimes_A N$ can be written, non-uniquly, as $\langle
								a,m_j\otimes n_j,a^{\prime}\rangle_n$. We thus get a bi-functor
								\begin{equation}
												 \_\otimes_A\_:\aset\times \aset\to \aset
												 \label{eq:2.13}
								\end{equation}
								giving a symmetric monoidal structure on $\aset$, with unit $A_{[1]}$.
								\vspace{.1cm}\\
								This symmetric monoidal structure on $\aset$ is closed, 
								\begin{equation}
												\aset(M\otimes_{A} N,K)\equiv \aset(M,\myhom_A(N,K)) 
												\label{eq:2.14}
								\end{equation}
								with the internal Hom functor  
								\begin{gather*}
								\myhom_A(\_,\_):(\aset)^{\op}\times\aset\to
								\aset \\
								\myhom_A(M,N):= \aset(M,N)
								\end{gather*}
								where the $A$-action on $\myhom_{A}(M,N)$ is given by
\[ \langle b,\varphi_j,d\rangle_n(m) := \langle b,\varphi_j(m),d\rangle_n, \quad b,d\in
A_n, \quad \varphi_j\in\myhom_A(M,N) \quad m\in M. \]
This in itself is a map of $\aset$ because of commutativity, and all the other properties
(\textsection (2.1): associativity, unit, naturality, commutativity) follow from their
validity in $N$. \\
The tensor product commute with extension of scalars: for $\varphi\in\cgr(B,A)$, and for
$M,N\in \bset$, we have 
\begin{equation*}
				\begin{aligned}[rl] 
				\varphi^*(M\otimes_B N) &\cong \varphi^*M\otimes_A\varphi^* N \\
				\varphi^* B_{[1]} &\cong A_{[1]}
				\end{aligned}
\end{equation*}
We have therefore also the adjunction formula
\begin{equation}
				\varphi_{*}\myhom_A(\varphi^{*}M,N) \cong \myhom_B(M,\varphi_* N)
				\label{eq:2.15}
\end{equation}
The tensor product is distributive over sums, 
\begin{equation}
				M\otimes_A\left( \coprod\limits_{i}N_i \right) \cong \coprod\limits_{i}\left( M\otimes_A N_i \right).
				\label{eq:2.16}
\end{equation}
and more generally commutes with colimits.
\section{Simplicial $\aset$, $\aset^{{{\Bnabla}}}$}
We denote by $\BDelta$ the simplicial category of finite ordered sets \linebreak $[n]=\left\{
0<1<\dots <n \right\}$, and monotone maps, and we let $\Bnabla=\BDelta^{\op}$ denote
the opposite category. We denote by $\aset^{\Bnabla}\equiv
(\aset)^{\BDelta^{\op}}$ the category of simplicial objects in $\aset$. It
has objects $M=(M_n)_{n\ge 0}\in \left( \set_0\right)^{\BDelta^{\op}}$, pointed
simplicial sets, with an $A$-action in each dimension $n\ge 0$, 
\begin{equation}
				\begin{array}[h]{ll}
				A_m\times \left( M_n \right)^m\times A_m \to M_n \\
				(a,m_i,a^{\prime}) \mapsto \langle a,m_i,a^{\prime}\rangle_m
				\end{array}
				\label{eq:3.1}
\end{equation}
compatible with the simplicial operations: 
\begin{equation}
				\begin{array}[H]{ll}
								\text{for } \ell\in\BDelta(n,n^{\prime}), \qquad & \ell^*:M_{n^{\prime}}\to M_n \;\text{ satisfies}
				\\\\
				& \ell^*(\langle a,m_i,a^{\prime}\rangle_m)=\langle a,\ell^*(m_i),a^{\prime}\rangle_m
				\end{array}
				\label{eq:3.2}
\end{equation}
The category $\asetnab$ is complete and co-complete. 
\begin{remark}
				 More generally, we have the category of
								simplicial-commutative-generalized-rings $\cgr^{\Bnabla}$, and for
								$A=\left( A^n \right)_{n\ge 0}\in \cgr^{\Bnabla}$ we have the category of
								simplicial $A\text{\normalfont-sets}$ $\asetnab$, with objects the pointed
								simplicial sets $M=\left( M_n \right)_{n\ge 0}$,  with compatible
								$A^n$-action in dimension $n\ge 0$, i.e. for $m\in \mathcal{F}$, 
								\begin{gather*}
												A_m^n\times \left( M_n \right)^m\times A_m^n\to M_n \\
												(a,m_i,a^{\prime})\mapsto \langle a, m_i, a^{\prime}\rangle_m^n
								\end{gather*}
								and for $\ell\in \BDelta(n,n^{\prime})$, 
								\[ \ell^{*}\left( \langle a,m_i,a^{\prime}\rangle_{m}^{n^{\prime}} \right)= 
												\langle
												\ell^{*}(a),\ell^{*}(m_i),\ell^{*}(a^{\prime})\rangle_{m}^{n}. 
								\]
\end{remark}
								\noindent While $\cgr^{\Bnabla}$ are important for
								''derived-arithmetical-geometry``, we shall concentrate here on $\cgr$, to
								keep the notation simpler. \par\vspace{.2cm}
								The category $\asetnab$ inherits a closed
								symmetric monoidal structure, 
								\begin{equation}
												\begin{array}[H]{ll}
												\_\otimes_A^\Bnabla\_: \asetnab\times\asetnab\to\asetnab \\\\
												\left( M_{.}\otimes_A^{\Bnabla}N_{.} \right)_n := M_n\otimes_A N_n
												
												\end{array}
												\label{eq:3.4}
								\end{equation}
		 \begin{equation}
						 \begin{array}[H]{ll}
										 \myhom_{A}^{\Bnabla}\left( \_,\_ \right): \left( \asetnab
										 \right)^{\op} \times \asetnab \to \asetnab \\\\
										 \myhom_A^{\Bnabla}(M_{.},N_{.})_n:=
										 \asetnab(M_{.}\otimes\phi_{A}^{*}(\Delta(n)_{+}),N_{.})
						 \end{array}
						 \label{eq:3.5}
		 \end{equation}
		 In particular, the category $\asetnab$ is tensored, co-tensored and enriched over
		 pointed simplicial sets. Denoting by $\phi_{A}\in \cgr(\mathbb{F},A)$ the unique map,
		 we have for $M_{.}\in \asetnab$, $K_{.}\in\text{$\mathbb{F}$-Set}^\Bnabla \equiv
		 \left( \set_0 \right)^{\BDelta^{\op}}$, 
		 \begin{equation}
						 \begin{array}[H]{lll}
										 (\romannumeral 1)\;\; M_{.}\otimes K_{.}:=
										 M_{.}\otimes_{A}^{\Bnabla}\phi_{A}^{*}K_{.} \\\\
										 (\romannumeral 2)\;\; M_{.}^{K_{.}} :=
										 \myhom_{A}^{\Bnabla}\left( \phi_{A}^{*}K_{.},M_{.} \right)\\\\
										 (\romannumeral 3)\;\; \text{Map}_{A}^{\Bnabla}\left( N_{.},M_{.} \right)
										 := \phi_{A^{*}}\myhom_{A}^{\Bnabla}\left( N_{.},M_{.} \right)\in
										 \set_0^{\Bnabla}
						 \end{array}
						 \label{eq:3.6}
		 \end{equation}
		 The category $\asetnab$ has a simplicial, cellular, Quillen model structure given by
		 cf.\cite{MR0223432}.
		 \begin{equation}
						 \begin{array}[h]{lll}
										 (\romannumeral 1) \text{ \underline{Fibrations}: } &
										 \mathcal{F}_{A} \equiv \phi_{*}^{-1}(\mathcal{F}_{\mathbb{F}}), &
										 \mathcal{F}_{\mathbb{F}}\equiv\text{Kan fibrations.}  \\\\
										 (\romannumeral 2) \text{ \underline{Weak-equivalence}: } &
										 \mathcal{W}_{A} \equiv
										 \phi_{*}^{-1}(\mathcal{W}_{\mathbb{F}}), &
										 \mathcal{W}_{\mathbb{F}}\equiv 
										 \hspace{-0.1cm} \left[ \hspace{-0.1cm}
										 \begin{array}[H]{l}
														 \parbox{2cm}{Weak-equivalence of (pointed) simplicial sets}
						 \end{array}\hspace{-0.25cm}\right] \\\\
		 
		 (\romannumeral 3) \text{ \underline{Cofibrations}: } & \mathcal{C}_A\equiv
										 \mathscr{L}(\mathcal{W}_{A}\cap \mathcal{F}_{A}), & 
										 \parbox[t][][t]{2.9cm}{the maps satisfying the left lifting property
														 with respect to the trivial fibrations.}
						 \end{array}
						 \label{eq:3.7}
		 \end{equation}
		 The cofibrations can also be characterized as the retracts of the \underline{free
		 maps}, where a map $f:M_{.}\to N_{.}$ is \underline{free} if there exists subsets
		 $V_n\subseteq N_n$ with $\ell^{*}(V_{n^{\prime}})\subseteq V_n$ for all surjective
		 $\ell\in\BDelta(n,n^{\prime})$, and such that $f_n$ induces isomorphism 
		 \begin{equation}
						 M_n\coprod \phi_{A}^{*}(V_n)\xrightarrow{\sim}N_n
						 \label{eq:3.8}
		 \end{equation}
		 where $\phi_{A}^{*}(V_n)$ is the free $\aset$ on $V_n$. This model structure is
		 cofibrantly generated 
		 \begin{equation}
						 \begin{array}[H]{ll}
										 (\romannumeral  1) & \mathcal{F}_A\equiv
										 \mathcal{R}\left\{\phi_{A}^{*}(\Lambda_{n+}^k)\to
										 \phi_{A}^{*}\big(\Delta(n)_{+}\big) \right\}_{0\le k\le n >0} \\\\
						 (\romannumeral  2) & \mathcal{W}_{A}\cap\mathcal{F}_A\equiv
						 \mathcal{R}\left\{\phi_{A}^{*}(\partial\Delta(n)_{+})\to\phi_A^{*}\left(
														 \Delta(n)_{+}
		 \right)\right\}_{n>0}
						 \end{array}
						 \label{eq:3.9}
		 \end{equation}
		 Where $\mathcal{R}\left\{ \_ \right\}$ are the map satisfying the Right lifting property
		 w.r.t.$\left\{ \_ \right\}$. \par
		 The model structure is compatible with the symmetric monoidal structure: \\
		 For  $\left\{ i_{.}\; :\; N_{.}\to N_{.}^{\prime} \right\}, \quad \left\{
		 j_{.}\; : \; M_{.}\to M_{.}^{\prime}\right\}$ in $\mathcal{C}_A$, we have

		 \begin{equation}
						 i_{.}\square j_{.}: (N_{.}\otimes_{A}^{\Bnabla}M_{.}^{\prime})
						 \coprod_{N_{.}\otimes_{A}^\Bnabla M_{.}} \left(
						 N_{.}^\prime\otimes_{A}^{\Bnabla} M_{.} \right)\longrightarrow
						 N_{.}^{\prime}\otimes_{A}^{\Bnabla}
						 M_{.}^{\prime}
						 \label{eq:3.10}
		 \end{equation}
		 is also in $\mathcal{C}_{A}$, and moreover, if $i_{.}$ is in $\mathcal{W}_A$, also
		 $i_{.}\square j_{.}$ is in $\mathcal{W}_{A}$. To see this compatibility of the monoidal
		 and model structures, we may assume $i_{.}$ and $j_{.}$ are free 
		 \[ N_{n}\coprod \phi_{A}^{*}(V_n)\xrightarrow{\sim} N_{n}^{\prime}, \quad M_n\coprod
		 \phi_{A}^{*}(W_{n})\xrightarrow{\sim} M_n^{\prime}\]
		 and than
		 \begin{equation}
						 \left( N_n\otimes_{A} M_n^{\prime} \coprod\limits_{N_n\otimes_A M_n}
										 N_n^{\prime}\otimes_A M_n
						 \right)\coprod \phi_{A}^{*}(V_n)\otimes_A
						 \phi_{A}^{*}(W_n)\xrightarrow{\sim} N_n^{\prime}\otimes_{A}M_{n}^{\prime}
						 \label{eq:3.11}
		 \end{equation}
		 and since $\phi_{A}^{*}(V_n)\otimes_{A}\phi_{A}^{*}(W_n)\equiv \phi_{A}^{*}(V_n\prod
		 W_n)$, we see that $i_{.}\square j_{.}$ is also free. \vspace{.2cm}\\
		 A homomorphism $\varphi\in\cgr(B,A)$ induces a Quillen adjunction,
		 \begin{equation}
						\begin{array}[H]{l}
								\begin{tikzpicture}
												\node (a) at (4,4) {$\asetnab$};
												\node (b) at (4,2) {$\bset^{\Bnabla}$};
												\draw [->] (a) to [out=330,in=30,right] node{$\varphi_{*}$} (b);
												\draw [<-] (a) to [out=210,in=150,left] node{$\varphi^{*}$} (b);
								\end{tikzpicture} 
						\end{array}
						\label{eq:3.12}
		 \end{equation}
		 When $A$ is a commutative ring, we have by the Dold-Kan correspondence $\asetnab\cong
		 \ch_{\ge 0}(A\mymod)$, and the model structure on $\asetnab$ corresponds to
		 the projective model structure on $\ch_{\ge 0}(A\mymod)$, which embeds in the
		 stable model structure $\ch_{+}(A\mymod)$. For general $A$, the model
		 structure on $\asetnab$ is not stable, and we shall stabilize it, (preserving the
		 symmetric monoidal structure), using (symmetric) spectra.

		 \section{Symmetric Spectra: $S_A^{\cdot}\mymod$ (cf.\cite{MR1860878})}
Let $\Sigma_n$ denote the symmetric group on $n$ letters, and let
$\Sigma_{.}=\coprod\limits_{n\ge 0}\Sigma_n$ denote the category of finite bijection. The
category $\Sigma_{.}$ is equivalent to the category 
\begin{equation}
				\text{Iso} (\mathcal{F})\equiv \text{Iso}(\fin_0)\equiv \text{Iso}(\fin)
				\label{eq:4.1}
\end{equation}
The category of symmetric sequences in $\asetnab$ is 
\begin{equation}
				\Sigma(A)\equiv \left( \asetnab \right)^{\Sigma_{.}} \equiv \left( \aset
				\right)^{\Sigma_{.}\times \Bnabla}\cong \left( \aset
				\right)^{\text{Iso}(\mathbb{F})\times \BDelta^{\text{op}}}
				\label{eq:4.2}
\end{equation}
It has objects $M^{\cdot}=\left\{ M^{n} \right\}_{n\ge 0}$, with $M^{n}\in  \left(\asetnab
\right)^{\Sigma_n}$ a simplicial- $\aset$ with an action of $\Sigma_{n}$. 
The category $\Sigma\left( A \right)$ is complete and co-complete.\vspace{.1cm}\\ 
The category $\Sigma\left( A \right)$ has a closed symmetric monoidal structure 
\begin{equation}
				\begin{array}[H]{ll}
				\_\otimes_{\Sigma\left( A \right)}\_ : \Sigma\left( A \right)\times
				\Sigma(A)\to \Sigma(A) \\\\
				\left( M^{\bullet}\otimes_{\Sigma(A)}N^{\bullet} \right)^{n}:=
				\coprod\limits_{p+q=n}\Sigma_{n}\underset{\Sigma_{p}\times\Sigma_{q}}{\times} \left(
								M^{p}\otimes_{A}^{\Bnabla} N^{q}
				\right)
				\end{array}
				\label{eq:4.3}
\end{equation}
Here the induction functor is the left adjoint of the forgetfull functor 
\[ \left( \aset \right)^{\Bnabla\times \Sigma_{n}}\to \left( \asetnab
\right)^{\Sigma_p\times \Sigma_q}, \]
and is given by 
\[ \Sigma_{n}\sigmapqtimes (M) :=
\coprod_{\Sigma_n/\Sigma_p\times\Sigma_q} M \]
Equivalently, writing $M^{\cdot},N^{\cdot}\in \Sigma(A)$ as functor $\text{Iso}(\fin)\to \asetnab$
we have 
\begin{equation}
				\left( M^{\cdot}\otimes_{\Sigma(A)}N^{\cdot} \right)^{n}:= \coprod_{n=n_0\coprod
				n_1} M^{n_0}\otimes_{A}^{\Bnabla} N^{n_1}
				\label{eq:4.4}
\end{equation}
the sum over all decomposition of $n$ as a disjoint union of subsets $n_0,n_1\subseteq n$.
The unit of this monoidal structure is the symmetric sequence 
\begin{equation}
				\characteristic_{A}:= \left( A_{[1]},0,0,0,\ldots \right)
				\label{eq:4.5}
\end{equation}
\begin{remark} \label{rem:4.6}
				Note that this monoidal structure is symmetric,
\begin{equation*}
				\mathfrak{T}_{M^{\cdot},N^{\cdot}}: M^{\cdot}\otimes_{\Sigma(A)}N^{\cdot}\cong
				N^{\cdot}\otimes_{\Sigma(A)}M^{\cdot}
\end{equation*}
This symmetry is clear in the formula (\ref{eq:4.4}),
$M^{n_o}\otimes_{A}^{\Bnabla}N^{n_1}\cong N^{n_1}\otimes_{A}^{\Bnabla}M^{n_0}$ but in
formula (\ref{eq:4.3}), the symmetry isomorphisms 
\[ \Sigma_{n} \sigmapqtimes \left( M^{p}\otimes_{A}^{\Bnabla}
				N^{q} \right) \cong \Sigma_{n}\sigmaqpt\left(
								N^{q}\otimes_{A}^{\Bnabla}M^p
\right) \]
involves the $(p,q)$-shuffle $\omega_{p,q}\in \Sigma_n$ that conjugates
$\Sigma_p\times\Sigma_q$ to $\Sigma_q\times\Sigma_p$. \vspace{.1cm}\\
\end{remark} 
\noindent The internal $\myhom$ is given by
\begin{equation}
				\begin{array}[H]{c}
								\myhom_{\Sigma(A)}(\_,\_): \Sigma(A)^{\op}\times
								\Sigma(A)\to \Sigma(A) \\\\
								\myhom_{\Sigma(A)}(M^{\cdot},N^{\cdot})^{n}:=\prod_{k\ge
								0}\myhom_{A}^{\Bnabla}(M^k,N^{k+n}) \\\\
								\Sigma(A)(M^{\cdot}\otimes_{\Sigma(A)}N^{\cdot},K^{\cdot})\cong
								\Sigma(A)\left( M^{\cdot},\myhom_{\Sigma(A)}(N^{\cdot},K^{\cdot}) \right)
				\end{array}
				\label{eq:4.7}
\end{equation}
For a homomorphism $\varphi\in\cgr(B,A)$, we have adjunction 
\begin{equation}
						\begin{array}[H]{l}
								\begin{tikzpicture}
												\node (SigmaA) at (4,4) {$\Sigma(A)$};
												\node (SigmaB) at (4,2) {$\Sigma(B)$};
												\draw [->] (SigmaA) to [out=330,in=30,right] node{$\varphi_{*}$}
												(SigmaB);
												\draw [<-] (SigmaA) to [out=210,in=150,left] node{$\varphi^{*}$}
												(SigmaB);
								\end{tikzpicture} 
						\end{array}
						\label{eq:4.8}
\end{equation}
and $\varphi^{*}$ is strict-monoidal 
\begin{equation}
				\varphi^{*}\left( M^{\cdot}\otimes_{\Sigma(B)}N^{\cdot} \right)\cong
				\varphi^{*}(M^{\cdot})\otimes_{\Sigma(A)}\varphi^{*}(N^{\cdot}), \quad
				\varphi^{*}\left( \characteristic_{B} \right)\cong \characteristic_{A}.
				\label{eq:4.9}
\end{equation}
The categogry $\Sigma(\mathbb{F})\equiv (\set_0)^{\Sigma_{.}\times \BDelta^{\op}}$ is the
usual category of symmetric sequence of pointed simplicial set, and in particular contains
the \underline{sphere-spectrum}: 
\begin{equation}
				S^{\cdot}_{\mathbb{F}} := \left\{ S^{n}=\underbrace{S^{1}\wedge\dots \wedge
				S^{1}}_{n} \right\}_{n\ge 0} 
				\label{eq:4.10}
\end{equation}
with the permutation action of $\Sigma_{n}$ on $S^{n}$. The sphere-spectrum is a monoid
object of $\Sigma(\mathbb{F})$, with multiplication
\begin{equation}
				\begin{array}[H]{c}
								m:
								S^{\cdot}_{\mathbb{F}}\otimes_{\Sigma(\mathbb{F})}S^{\cdot}_{\mathbb{F}}\to
								S^{\cdot}_{\mathbb{F}} \\\\
								m(S^{n}\otimes_{\mathbb{F}}^{\Bnabla}S^{m}) \equiv m(S^n\wedge S^m)
								\equiv S^{n+m} 
				\end{array}
				\label{eq:4.11}
\end{equation}
Note that it is a \underline{commutative} monoid,
$m=m\circ\mathfrak{T}_{{S_{\mathbb{F}}^{\cdot}},S^{\cdot}_{\mathbb{F}}}$, cf. remark
(\ref{rem:4.6}). \\
The unit is given by the embedding 
\begin{equation}
				\varepsilon: \characteristic_{\mathbb{F}}\equiv \left(
				\mathbb{F}_{[\unit]},0,0,\dots \right) \equiv \left( S^{0},0,0,\dots
				\right)\hookrightarrow S^{\cdot}_{\mathbb{F}}
				\label{eq:4.12}
\end{equation}
We write $S^{\cdot}_{A}=\phi^{*}_{A}S^{\cdot}_{\mathbb{F}}$ for the corresponding
commutative monoid object of $\Sigma(A)$. We let $S^{\cdot}_{A}\mymod\subseteq
\Sigma(A)$ denote the sub-category of $S_{A}^{\cdot}$-modules, this is the category of
``symmetric spectra''. It has objects the symmetric sequences \linebreak $M^{\cdot}=\left\{
				M^{n}
\right\}\in \Sigma(A)$, thogether with associative unital 
 $S_{A}^{\cdot}$-action \linebreak
$S_{A}^{\cdot}\otimes_{\Sigma(A)}M^{\cdot}\xrightarrow{m}M^{\cdot}$, or equivalently,
associative unital, $\Sigma_{p}\times \Sigma_{q}\hookrightarrow \Sigma_{p+q}$ covariant,
action $S^{p}\wedge M^{q}\to M^{p+q}$. The maps in $S_{A}^{\cdot}\mymod$ are the maps in
$\Sigma(A)$ that preserve the $S_{A}^{\cdot}$-action. The category $S_{A}^{\cdot}\mymod$ is
complete and co-complete. \vspace{.1cm}\\
The category $S_{A}^{\cdot}\mymod$ has a closed symmetric monoidal structure 
\begin{equation}
				\begin{array}[H]{c}
								\_\otimes_{S^{\cdot}_{A}}\_: S^{\cdot}_{A}\mymod\times
								S_{A}^{\cdot}\mymod\to S^{\cdot}_{A}\mymod \\\\
								M^{\cdot}\otimes_{S_{A}^{\cdot}}N^{\cdot} := \text{Cok}\left\{
												M^{\cdot}\otimes_{\Sigma(A)}S_{A}^{\cdot}\otimes_{\Sigma(A)}N^{\cdot}\overset{\xrightarrow{m\otimes\text{id}_{N^{\cdot}}}}{
								\xrightarrow[\text{id}_{M^{\cdot}}\otimes m]{}} M^{\cdot}\otimes_{\Sigma(A)}N^{\cdot} \right\}
								\end{array}
								\label{eq:4.13}
\end{equation}
The unit is
\begin{equation}
				S_{A}^{\cdot} := \left\{ \phi_{A}^{*}S^{0},\phi_{A}^{*}S^{1},\ldots ,
				\phi_{A}^{*}S^{n},\ldots \right\}. 
				\label{eq:4.14}
\end{equation}
The internal $\myhom$ is given by 
\begin{equation}
				\begin{array}[H]{cc}
				\myhom_{S_{A}}(\_,\_): \left( S_{A}^{\cdot}\mymod \right)^{\op}\times
				S_{A}^{\cdot}\mymod\to S_{A}^{\cdot}\mymod \\\\
				\myhom_{S_A}(M^{\cdot},N^{\cdot}):= \kernel \left\{
				\myhom_{\Sigma(A)}(M^{\cdot},N^{\cdot})\rightrightarrows 
				\myhom_{\Sigma(A)}\left( S^{\cdot}_{A}\underset{{\Sigma(A)}}{\otimes} M^{\cdot},N^{\cdot}
\right)\right\}.
				\end{array}
				\label{eq:4.15}
\end{equation}
\begin{equation}
				S_A\mymod\left(M^{\cdot}\otimes_{S^{\cdot}_{A} }N^{\cdot},K^{\cdot}  \right)
				\equiv S_{A}\mymod\left(
				M^{\cdot},\myhom_{S^{\cdot}_{A}}(N^{\cdot},K^{\cdot}) \right)
				\label{eq:4.16}
\end{equation}
The category $S_{A}^{\cdot}\mymod$ is tensored, co-tensored, and enriched over pointed
simplicial sets $\mathbb{F}\text{-Set}^{\Bnabla}:$ For
$K\in\mathbb{F}\text{-Set}^{\Bnabla}$, $M^{\cdot}\in S_{A}\mymod$, we have
$M^{\cdot}\otimes K$, $\left( M^{\cdot} \right)^{K}\in S^{\cdot}_{A}\mymod$ where 
\begin{equation}
				\left( M\otimes K \right)^{n} = M^{n}\otimes_{\mathbb{F}}^{\Bnabla}K =
				M^{n}\wedge K
				\label{eq:4.17}
\end{equation}
\begin{equation}
				\left( M^{\cdot^{K}} \right)^{n} = \left( M^{n} \right)^{K} \equiv
				\mathbb{F}\text{-Set}^{\Bnabla}\left( K\wedge \Delta(\_)_{+},M^{n} \right)
				\label{eq:4.18}
\end{equation}
The enrichment is given via the mapping space 
\begin{equation}
				\text{Map}_{S^{\cdot}_{A}}\left( M^{\cdot},N^{\cdot} \right)_{n} \equiv
				S_{A}^{\cdot}\mymod\left( M^{\cdot}\otimes \Delta(n)_{+},N^{\cdot} \right)\in
				\mathbb{F}\text{-Set}^{\Bnabla} \equiv \left( \text{\set}_0
				\right)^{\BDelta^{\op}}
				\label{eq:4.19}
\end{equation}
We have the adjunctions, 
\begin{equation}
				\begin{array}{c}
								\begin{array}[H]{lll}
												S^{\cdot}_{A}\mymod\left( M^{\cdot}\otimes K, N^{\cdot} \right) &\equiv&
												S_{A}^{\cdot}\mymod\left( M^{\cdot},\left( N^{\cdot} \right)^{K} \right) \\
												& \equiv&\mathbb{F}\text{-Set}^{\Bnabla}\left(
												K,\text{Map}_{S^{\cdot}_{A}}(M^{\cdot},N^{\cdot}) \right) 
								\end{array} \\\\
								\begin{array}[h]{c}
												\begin{tikzpicture}
																\node (samod) at (4,4) {$S_{A}^{\cdot}\mymod$};
																\node (SigmaB) at (4,2) {$\mathbb{F}\text{-set}^{\Bnabla}$};
																\draw [->] (samod) to [out=330,in=30,right]
																node{$\text{Map}_{S^{\cdot}_{A}}\left( M^{\cdot},\_ \right)$}
																(SigmaB);
																\draw [<-] (samod) to [out=210,in=150,left]
																node{$M^{\cdot}\otimes_{\_}$}
																(SigmaB);
												\end{tikzpicture} 
								\end{array}
				\end{array}
						\label{eq:4.20}
\end{equation}
\begin{equation}
				\begin{array}{lll}
								\text{Map}_{S^{\cdot}_{A}	}\left( M^{\cdot}\otimes K, N^{\cdot} \right)
								&\equiv&
								\text{Map}_{S^{\cdot}_{A}}\left( M^{\cdot},\left( N^{\cdot} \right)^{K}
								\right) \\
								&\equiv& \text{Map}_{S^{\cdot}_{A}}\left( M^{\cdot},N^{\cdot} \right)^{K}
				\end{array}
				\label{eq:4.21}
\end{equation}
Taking $K=\text{Map}_{S^{\cdot}_{A}}(M^{\cdot},N^{\cdot})$, the
$\text{id}_{\text{Map}(M^{\cdot},N^{\cdot})}\in
\text{Map}_{S^{\cdot}_{A}}(M^{\cdot},N^{\cdot})^{K}$, on r.h.s. (\ref{eq:4.21}) corresponds
to the evaluation map on $l.h.s$ (\ref{eq:4.21}) \vspace{.1cm}
\begin{equation}
				\text{ev}_{M^{\cdot},N^{\cdot}}:
				M^{\cdot}\otimes_{\mathbb{F}}\text{Map}_{S^{\cdot}_{A}}(M^{\cdot},N^{\cdot})\to
				N^{\cdot} 
				\label{eq:4.22}
\end{equation}
Taking $K=\text{Map}\left( M^{\cdot},L^{\cdot} \right)\otimes_{\mathbb{F}}\text{Map}\left(
L^{\cdot},N^{\cdot} \right)$, the map
$$ \text{ev}_{L^{\cdot},N^{\cdot}}\circ
\left(\text{ev}_{M^{\cdot},L^{\cdot}}\otimes\text{id}_{\text{Map}(L^{\cdot},N^{\cdot})}\right)$$
on the $\ell$.h.s of (\ref{eq:4.21}) corresponds to the composition map, on r.h.s. 
\begin{equation}
				\text{comp}_{L^{\cdot}}:
				\text{Map}(M^{\cdot},L^{\cdot}) \otimes_{\mathbb{F}}
				\text{Map}(L^{\cdot},N^{\cdot})\to \text{Map}(M^{\cdot},N^{\cdot})
				\label{eq:4.23}
\end{equation}
which is associative and unital. \\
We have the embedding $\left( A\text{-Set}^{\Bnabla} \right)^{\Sigma_{n}}\hookrightarrow
\Sigma(A)$, $M\mapsto M[n]$, and the following adjunctions 
\begin{equation}
				\begin{array}[H]{l}
								\begin{tikzpicture}
												\node (B0) at (4,7) {$S_{A}^{\cdot}\mymod$};
												\node (B1) at (4,4) {$\left( \asetnab \right)^{\Sigma_{n}}$};
												\node (B2) at (4,1) {$\asetnab$};
												\node (B00) at (4.4,6.6) {$M^{\cdot}$};
												\node (B01) at (4.4,4.8) {$M^{n}$};
												\draw [->] (B1) to [out=30,in=-10,right] (B0);
												\draw [->] (B1) to [out=150,in=195,left] (B0);
												\draw [<-] (B1) to [out=330,in=30,right] (B2);
												\draw [<-] (B1) to [out=210,in=150,left] (B2);
												\draw [->] (B1) to [right] 
																node{$U$} (B2);
												\draw [->] (B0) to [left] 
												node{$Ev^{n}$} (B1);
												\draw [Bar->] ($(B00)-(0.1,0.2)$) to ($(B01)-(0.1,-0.2)$);
												\node (A0) at (0.5,7) {$S_{A}^{\cdot}\otimes_{\Sigma(A)}M[n]$};
												\node (A01) at (0.5,4.8) {$M$};
												\node (A02) at (0.5,3.5) {$\Sigma_{n}\times M=\coprod\limits_{\Sigma_{n}}M$};
												\node (A03) at (0.5,1) {$M$};
												\draw [Bar->] (A01) to (A0);
												\draw [Bar->] (A03) to (A02);

												\node (C0) at (7.5,7) {$\myhom_{\Sigma(A)}(S_{A}^{\cdot},M[n])$};
												\node (C01) at (7.5,4.8) {$M$};
												\node (C02) at (7.5,3.5) {$M^{\Sigma_{n}}=\prod\limits_{\Sigma_n}M$};
												\node (C03) at (7.5,1) {$M$};
												\draw [Bar->] (C01) to (C0);
												\draw [Bar->] (C03) to (C02);
								\end{tikzpicture} 
								\end{array}
				\label{eq:4.24}
\end{equation}
The {\underline{``free-$S^{\cdot}_{A}$-module of level-$n$}}'' on $M\in \asetnab$ is the
composition of the left adjoints of (\ref{eq:4.24}), 
\begin{equation}
				F_{n}\left( M \right):=S_{A}^{\cdot}\otimes_{\Sigma(A)}\left( \Sigma_{n}\times M
				\right)[n] = \left( 0,\dots, 0, \Sigma_n\times
				M,\dots,\Sigma_{n+p}\underset{\Sigma_{p}}{\times}
				S^{p}\otimes M ,\dots \right).
				\label{eq:4.25}
\end{equation}
We have the projective model structure on $S_{A}^{\cdot}\mymod$, compatible with the
symmetric monoidal structure, 
\begin{equation}
				\begin{array}[H]{lll}
								\text{\underline{Fibrations}:} & \mathcal{F}_{S^{\cdot}_A}^{\text{lev}} =  \left\{
												f^{\cdot}\in S^{\cdot}_{A}\mymod(M^{\cdot},N^{\cdot}),f^{n}\in\mathcal{F}_{A}\text{ all } n\ge 0 \right\} \\\\
				\text{\underline{Weak equivalences}:} & \mathcal{W}_{S^{\cdot}_{A}}^{\text{lev}}= \left\{
								f^{\cdot}\in S^{\cdot}_{A}\mymod (M^{\cdot},N^{\cdot}),f^{n}\in
								\mathcal{W}_{A}\text{ all }
n\ge 0 \right\}\\\\
\text{\underline{Cofibrations}:} & \mathcal{C}_{S^{.}_{A}} = \mathscr{L}\left\{
\mathcal{W}^{\text{lev}}_{S^{\cdot}_{A}}\cap \mathcal{F}_{S^{\cdot}_{A}}^{\text{lev}} \right\}
				\end{array}
				\label{eq:4.26}
\end{equation}
It is left proper, simplicial, cofibrantly generated by 
\begin{equation}
				\begin{array}[H]{lll}
				J &=& \coprod\limits_{\overset{m\ge 0}{n\ge 1}} F_m\left( \phi_{A}^{*}\partial
				\Delta(n)_{+}\hookrightarrow \phi_{A}^{*}\Delta(n)_{+} \right) \\\\
				I &=& \coprod\limits_{\overset{m\ge 0}{0\le k\le n\ge 1}} F_{m}\left(
				\phi_{A}^{*}\Lambda_{k_{+}}^{n}\hookrightarrow \phi_{A}^{*}\Delta(n)_{+} \right)
				\end{array}
				\label{eq:4.27}
\end{equation}
A map $f^{\cdot}\in S^{\cdot}_{A}\mymod\left( M^{\cdot},N^{\cdot} \right)$ is a
cofibrantion, $f\in \mathcal{C}_{S^{\cdot}_{A}}$, if and only if $Ev^{n}\left(
				f^{\cdot}\square j
\right): M^{n}\coprod_{\mathscr{L}^{n}M^{\cdot}}\mathscr{L}^{n}N^{\cdot}\to N^{\cdot}$ is
in $\mathcal{C}_{A}$ and $\Sigma_{n}$ acts freely away from its image, for all $n\ge 0$,
where the ''latching functor`` $\mathscr{L}^{n}$ is 
\begin{equation}
				\mathscr{L}^n M^{\cdot} = Ev^{n}\left( M\otimes_{\Sigma(A)} \left(
				S^{\cdot}_{A}/S^{0}_{A} \right) \right) = \coprod\limits_{0\le k < n}
				\Sigma^{n}\underset{\Sigma_{k}\times\Sigma_{n-k}}{\times}  \left(
				M^{k}\otimes_{A}^{\Bnabla} S^{n-k}_{A} \right)
				\label{eq:4.28}
\end{equation}
We say $M^{\cdot}\in S^{\cdot}_{A}\mymod$ is \underline{an $\Omega$-spectrum}, or a
fibrant object $M^{\cdot}\in \left( S^{\cdot}_{A}\mymod \right)_{\mathcal{F}}$ if
$M^{\cdot}$ is levelwise fibrant, and if the adjoint of the action map  
\[ m^{1,n}:
S_{A}^{1} \otimes_{A}^{\Bnabla} M^{n}\to M^{n+1} \]
is a weak equivalence
\begin{equation}
				\left( m^{1,n} \right)^{\natural} : M^{n}\xrightarrow{\sim}\myhom_{A}^{\nabla}\left(
				S_{A}^{1},M^{n+1} \right) = (M^{n+1})^{S^{1}}:= \Omega M^{n+1}
				\label{eq:4.29}
\end{equation}
The \underline{stable} model structure on $S^{\cdot}_{A}\mymod$ is a Bousfield localization
of the projective model structure, having the same cofibrantions, but with the fibrant
objects being $\left( S^{\cdot}_{A}\mymod \right)_{\mathcal{F}}$ 
\begin{equation}
				\begin{array}[H]{ll}
								\text{\underline{Cofibrations}:}			&
								\mathcal{C}_{S^{\cdot}_{A}}= \mathscr{L}\left\{
												\mathcal{W}_{S^{\cdot}_{A}}^{\text{lev}}\cap
				\mathcal{F}_{S^{\cdot}_{A}}^{\text{lev}} \right\}												 \\\\
				\text{\underline{Weak equivalences}:} &	 \mathcal{W}_{S^{\cdot}_{A}}=\left\{ 
								\begin{array}[H]{lr}
												f^{\cdot}\in S^{\cdot}_{A}\mymod\left( M^{\cdot},N^{\cdot}
								\right), 
								\\ \text{Map}_{S^{\cdot}_{A}}(N^{\cdot},
								X^{\cdot})\xrightarrow{\sim}\text{Map}_{S^{\cdot}_{A}}(M^{\cdot},X^{\cdot})\in
								\mathcal{W}_{\mathbb{F}}  \\
								\text{\hspace{2cm} for all } X^{\cdot}\in \left(
												S^{\cdot}_{A}\mymod
								\right)_{\mathcal{F}} 
								\end{array}
				\right\} \\\\
								\text{\underline{Fibrations}:} &	\mathcal{F}_{S^{\cdot}_{A}} =  \mathcal{R}\left\{
												\mathcal{C}_{{S}^{\cdot}_{A}} \cap \mathcal{W}_{{S}^{\cdot}_{A}}
								\right\}
				\end{array}
				\label{eq:4.30}
\end{equation}
It is left proper, simplicial, cellular, Quillen Model structure compatible with its
symmetric monoidal structure. It is moreover \underline{stable}. We let
\begin{equation}
				\mathbb{D}(A)=\text{Ho}\left( S^{\cdot}_{A}\mymod \right)\cong
				S^{\cdot}_{A}\mymod\left[ \mathcal{W}_{S^{\cdot}_{A}}^{-1} \right]
				\label{eq:4.31}
\end{equation}
denote the associated homotopy category, this is the {\underline{derived}} category of
$\aset$. It is a triangulated symmetric monidal category. \vspace{.1cm}\\ 
We have the following Quillen adjunctions
\begin{equation}
				\begin{array}[H]{l}
								\begin{tikzpicture}
												\node (A) at (4,4) {$S_{A}^{\cdot}\mymod$};
												\node (B) at (4,2){$S_{A}^{\cdot}\mymod$};
												\draw [->] (A) to [out=330,in=30,right]
												node{$\Omega$}
												(B);
												\draw [<-] (A) to [out=210,in=150,left]
												node{$\_\underset{\mathbb{F}}{\otimes}S^{1}$}
												(B);

												\node (A0) at (9,4) {$S_{A}^{\cdot}\mymod$};
												\node (B0) at (9,2){$S_{A}^{\cdot}\mymod$};
												\draw [->] (A0) to [out=330,in=30,right]
												node{$r$}
												(B0);
												\draw [<-] (A0) to [out=210,in=150,left]
												node{$\ell$}
												(B0);
								\end{tikzpicture} 
				\end{array}
				\label{eq:40.32}
\end{equation}
where 
\[
				\begin{array}[H]{lll}
								\left( r M^{\cdot} \right)^{n}:= M^{n+1}, & r M^{\cdot} :=
\myhom_{S^{\cdot}_{A}}(F_{1}S^{\circ}_{A},M^{\cdot}) \\\\
\left( \ell M^{\cdot} \right)^{n}:= \Sigma_{n}\underset{\Sigma_{n-1}}{\times} M^{n-1},
\qquad &
\ell M^{\cdot} := F_{1}S_{A}^{\circ}\otimes_{S^{\cdot}_{A}}M^{\cdot} \\\\
\Omega M^{\cdot}:=\left( M^{\cdot} \right)^{S^{1}} \\\\
M^{\cdot}\otimes_{\mathbb{F}}S^{1} = M^{\cdot}\otimes_{A}S^{1}_{A} 
\end{array} \]
They induce Quillen equivalence, and we have the inverse equivalences of
$\mathbb{D}(A)$, 
\begin{equation}
				\R \Omega\simeq \mathbb{L}\ell, \qquad
				\_\otimes^{\mathbb{L}}_{\mathbb{F}}S^{1}
				\simeq \R r
				\label{eq:4.33}
\end{equation}
\ \vspace{-.2cm}\\
Given a homomorphism $\varphi\in \cgr(B,A)$ we get an induced adjunction 
\begin{equation}
				\begin{array}[H]{l}
								\begin{tikzpicture}
												\node (A) at (4,4) {$\mathbb{D}(A) $};
												\node (B) at (4,2){$\mathbb{D}(B)$};
												\draw [->] (A) to [out=330,in=30,right]
												node{$\R \varphi_{*}$}
												(B);
												\draw [<-] (A) to [out=210,in=150,left]
												node{$\mathbb{L}\varphi^{*}$}
												(B);
								\end{tikzpicture} 
				\end{array}
				\label{eq:4.34}
\end{equation}
and $\mathbb{L}\varphi^{*}$ is a monoidal functor commuting with $\mathbb{L}\ell$ and
$\R r$. We have the adjunction formula for $M^{\cdot}\in\mathbb{D}(B)$, $N^{\cdot}\in
\mathbb{D}(A)$,
\begin{equation}
				\R \varphi_{*}\left( \R\myhom_{S^{\cdot}_{A}}\left(
				\mathbb{L}\varphi^{*}M^{\cdot},N^{\cdot} \right) \right)\equiv
				\R\myhom_{S^{\cdot}_{B}} \left( M^{\cdot},\R\varphi_{*}N^{\cdot} \right)
				\label{eq:4.35}
\end{equation}
\section{Global theory}
For $A\in \text{GR}$, commutative or not, the sub-$A$-sets of $A_{[1]}$ are the
{\underline{ideals}} of $A$. Such an ideal $\plusideal\subseteq A_{[1]}$ is called a
{\underline{$+$ideal}} if it is generated as an ideal by its subset of
{\underline{symmetric elements}} $\plusideal^{+} = \plusideal \cap A^{+}_{[1]}=\left\{
a\in \plusideal\; , \; a^t = a \right\}$. A $+$ideal
$\plusprimeideal\subseteq A_{[1]}$ is called {\underline{+prime}} if the set
$S_{\plusprimeideal}=A_{[1]}^{+}\setminus\plusprimeideal$ is multiplicatively closed. The
set of $+$primes of $A$, $\text{spec}^{+}(A)$, is a {\underline{compact}} (every open
cover has a finite subcover), Sober or Zariski (every closed irreducible subset has a unique generic point)
topological space with closed subsets
$V^{+}(\plusideal):=\left\{ \plusprimeideal \supseteq \plusideal \right\}$, $\plusideal$ a
$+$ideal, and basic open subset  $D^{+}(s):=\left\{ \plusprimeideal
\not\ni s \right\}$, $s=s^{t}\in A_{[1]}^{+}$. By
localization one obtains a sheaf $\bigo_{A}$ of $\text{GR}$ over $\text{spec}^{+}(A)$, with
stalks $\bigo_{A}|_{\plusprimeideal}\cong A_{\plusprimeideal} =
S_{\plusprimeideal}^{-1} A$, and with sections $\bigo_{A}\left( D^{+}(s)
\right)\cong A\left[ \frac{1}{s} \right]$, cf. \cite{MR3605614}. \vspace{.2cm}\\
A \underline{generalized scheme} $(X,\bigo_{X})\in \text{GSch}$, is a
topological space $X$, a sheaf of generalized rings $\bigo_{X}$ over $X$, with local
stalks, and we have an open cover $X=\bigcup\limits_{i\in I} U_i$, with
$(U_i,\bigo_{X}|_{U_i})\cong \text{spec}^{+}\left( \bigo_{X}(U_{i}) \right)$. We can take
$I=\text{Aff}/X$, the collection of all affine open subsets of $X$, and we view $\text{Aff}/X$ as a
category with respect to inclusions. We denote by $\phi_{V,U}:\bigo_{X}(U)\to
\bigo_{X}(V)$ the restriction homomorphism for $V\subseteq U$, $\phi_{W,V}\circ
\phi_{V,U}=\phi_{W,U}$, $\phi_{U,U}=\text{id}_{\bigo_{X}(U)}$. \\
We shall assume all $\bigo_{X}(U)$ are commutative (1.26). 
\begin{definition} \label{def:5.1}
				\begin{textit} A $\pre\bigo_{X}\text{\normalfont-set}$ (resp.
								$\pre\bigo_{X}\text{\normalfont-Set}^{\Bnabla}$, $\pre\Sigma(\bigo_{X})$,
				$\pre S^{.}_{\bigo_{X}}\mymod$) is an assignment for each
								open $U\subseteq X$ of an $\bigo_X(U)$-Set $M(U)$ (resp.
								$\bigo_{X}(U)\text{-Set}^{\Bnabla}$, $\Sigma\left( \bigo_{X}(U) \right)$,
				$S^{\cdot}_{\bigo_{X}(U)}\mymod$) and for $V\subset U$, an
				$\bigo_X(U)\text{\normalfont-map}$
				$r_{V,U}: M(U)\to M(V)$, where $M(V)$ is acted on by $\bigo_X(U)$ via
				$\phi_{V,U}$, and $r_{W,V}\circ r_{V,U}=r_{W,U}$,
				$r_{U,U}=\text{id}_{M(U)}$.
				\end{textit}
\end{definition}
				The maps are the natural transformations
				commuting with the $\bigo_X$-action (resp. and with the simplicial  $\Bnabla$;
								symmetric
				$\Sigma = \coprod\limits_{n}\Sigma_n$; sphere $S^{\cdot}=\left\{ S^{n} \right\}$; actions) so that we have categories: 
\begin{equation}
				\begin{array}[H]{lll}
				& \pre\bigo_{X}\text{-Set} \\\\
				\text{resp.} & \pre\bigo_{X}\text{-Set}^{\Bnabla} \equiv  \left(
								\pre\bigo_{X}\text{-Set} \right)^{\BDelta^{\text{op}}} \\\\
								\ & \pre S^{\cdot}_{\bigo_{X}}\mymod \subseteq
								\pre\Sigma\left( \bigo_{X} \right) \equiv
								\left( \pre\bigo_{X}\text{-Set}
								\right)^{\BDelta^{\text{op}}\times \coprod\limits_{n}\Sigma_n}
				\end{array}
				\label{eq:5.2}
\end{equation}
All these categories are complete and co-complete (limits and colimits are created
sectionwise over open subsets). They have a closed symmetric monoidal structure defined
sectionwise which we shall denote by $\tilde{\otimes}_{\bigo_{X}}$ (resp. $\tilde{\otimes
				}_{\bigo_{X}^{\Bnabla}}$ ; $\tilde{\otimes}_{\Sigma(\bigo_X)}$ ; $ 
\tilde{\otimes}_{S^{\bullet}_{\bigo_{X}}}$), 
\begin{equation}
				\begin{array}[H]{lll}
				(M\tilde{\otimes}_{\bigo_{X}}N)(U) := M(U)\otimes_{\bigo_X(U)} N(U)
				\end{array}
				\label{eq:5.3}
\end{equation}
The internal Hom will be denoted by $\myhom_{\bigo_X}$ (resp.
$\myhom_{\bigo_X}^{\Bnabla}$; $\myhom_{\Sigma(\bigo_X)}$
$\myhom_{S^{\cdot}_{\bigo_X}}$) and is defined as the equalizer 
\begin{equation}
				\begin{array}[H]{lll}
\myhom_{\bigo_{X}}\left( M,M^{\prime} \right)(U) & \xrightarrow{\quad}
				\prod\limits_{\mathcal{V}\subseteq U} \myhom_{\bigo_{X}(\mathcal{V})}\left(
				M(\mathcal{V}),M^{\prime}(\mathcal{V}))  \right) \\\\
				& \rightTwoArrows 
				\prod\limits_{\mathcal{W}\subseteq \mathcal{V}\subseteq U}
				\myhom_{\bigo_X(\mathcal{V})}\left(
				M(\mathcal{V}),M^{\prime}(\mathcal{W}) \right)  \\\\
				& \ \hspace{2.4cm} \left\{ \varphi_{\mathcal{V}} \right\} \mapstoto 
				\begin{array}[H]{ll}
								\varphi_{\mathcal{W}}\circ r_{\mathcal{W},\mathcal{V}} \\
								r^{\prime}_{\mathcal{W},\mathcal{V}}\circ\varphi_{\mathcal{V}}
				\end{array}
				\label{eq:5.4}
				\end{array}
\end{equation}
We have the full subcategories of sheaves 
\begin{equation}
				\begin{array}[H]{ll}
								\ & \bigo_{X}\text{-Set}\subseteq \pre\bigo_{X}\text{-Set} \\\\
								\text{resp.} & \bigo_{X}\text{-Set}^{\Bnabla}\subseteq
				\pre\bigo_{X}\text{-Set}^{\Bnabla} \\\\
				\ & \Sigma(\bigo_{X})\subseteq \pre\Sigma(\bigo_{X}) \\\\
				\ & S^{\cdot}_{\bigo_X}\mymod\subseteq \pre S^{\cdot}_{\bigo_{X}}\mymod
				\label{eq:5.5}
				\end{array}
\end{equation}
where a pre-sheaf $M(U)$ is a {\underline{sheaf}} if for all open covers
$U=\bigcup\limits_{i\in I}U_i$, 
\begin{equation}
				\begin{array}[H]{rll}
				M(U)\xrightarrow{\quad}\prod_{i} M(U_i)
				& \rightTwoArrows  \prod_{i_0,i_1}M\left( U_{i_0}\cap U_{i_1} \right)  \\\\
				\left\{ m_i \right\} & \mapstoto \ & \  \vspace{-0.61cm}\\
				&    \hspace{0.85cm}  \begin{array}[H]{l}
								r_{u_{i_0} \cap u_{i_1}, u_{i_0}}(m_{i_0}) \\
								r_{u_{i_0} \cap u_{i_1}, u_{i_1}}(m_{i_1}) 
				\end{array}
				\end{array}
				\label{eq:5.6}
\end{equation}
is an equalizer. The above embeddings of categories (\ref{eq:5.5}) have left-adjoints
$M\mapsto M^{\natural}$ the \underline{sheafification functor}. These subcategories of
sheaves are complete (limits are created at the pre-sheaf level), and cocomplete (colimits
are obtained via the pre-sheaf colimit followed by sheafification). \vspace{.1cm}\\
The categories of sheaves are closed symmetric monoidal via 
\begin{equation}
				\begin{array}[H]{rll}
								M\otimes_{\bigo_{X}}N &:= \left( M\tilde{\otimes}_{\bigo_X}N
								\right)^{\natural} & \text{(using sheafification)} \\\\
								\myhom_{\bigo_{X}}\left( M,N \right) &:= \myhom_{\bigo_{X}}\left(
								M,N \right)
				\end{array}
				\label{eq:5.7}
\end{equation}
Notice that the internal pre-sheaf Hom of two sheaves is also a sheaf. \vspace{.1cm}\\
Given a map of generalized schemes $f\in\text{GSch}\left( X,Y \right)$, we get adjunctions 
\begin{equation}
				\begin{array}[h]{l}
								\begin{tikzpicture}
												\node (A0) at (1,4) {$\bigo_{X}\text{-Set}$};
												\node (B0) at (1,2){$\bigo_{Y}\text{-Set}$};
												\draw [->] (A0) to [out=330,in=30,right]
												node{$f_{*}$}
												(B0);
												\draw [<-] (A0) to [out=210,in=150,left]
												node{$f^{*}$}
												(B0);

												\node (A) at (5,4) {$\bigo_{X}\text{-Set}^{\Bnabla}$};
												\node (B) at (5,2){$\bigo_{Y}\text{-Set}^{\Bnabla}$};
												\draw [->] (A) to [out=330,in=30,right]
												node{$f_{*}$}
												(B);
												\draw [<-] (A) to [out=210,in=150,left]
												node{$f^{*}$}
												(B);

								\end{tikzpicture}  \\\\
								\begin{tikzpicture}
												\node (A2) at (1,4) {$\Sigma (\bigo_{X})$};
												\node (B2) at (1,2){$\Sigma (\bigo_{Y})$};
												\draw [->] (A2) to [out=330,in=30,right]
												node{$f_{*}$}
												(B2);
												\draw [<-] (A2) to [out=210,in=150,left]
												node{$f^{*}$}
												(B2);

												\node (A3) at (5,4) {$S_{\bigo_{X}}^{\cdot}\mymod$};
												\node (B3) at (5,2){$S_{\bigo_{Y}}^{\cdot}\mymod$};
												\draw [->] (A3) to [out=330,in=30,right]
												node{$f_{*}$}
												(B3);
												\draw [<-] (A3) to [out=210,in=150,left]
												node{$f^{*}$}
												(B3);
								\end{tikzpicture}
				\end{array} 
				\label{eq:5.8}
\end{equation}
satisfying the usual properties, such as the existence of canonical equivalences
\begin{equation}
				\begin{array}[H]{c}
								\left( g\circ f \right)_{*}\xrightarrow{\;\sim\;} g_{*}\circ f_{*} \qquad
								\left( g\circ f \right)^{*}\xrightarrow{\;\sim\;} f^{*}\circ g^{*} \\\\
								f^{*}\left( M\otimes_{\bigo_{Y}} N \right) \xrightarrow{\;\sim\;}
								f^{*} M\otimes_{\bigo_{X}} f^{*}N \\\\
								f_{*}\myhom_{\bigo_{X}}\left( f^{*}M,N \right)\xrightarrow{\;\sim\;}
								\myhom_{\bigo_{Y}}\left( M,f_{*}N \right). 
				\end{array}
\label{eq:5.9}
\end{equation}
Given $A\in \cgr$, $S\subseteq A^{+}_{[1]}$ a multiplicative subset, and $M\in \aset$, we
can localize $M$ to obtain $S^{-1}M\in S^{-1}A\text{-Set}$ in the usual way: $S^{-1}M :=
\left( M\times S \right)/\simeq$, with the equivalence relation $\simeq$ given by 
\begin{equation}
				\left( m_1,s_1 \right)\simeq \left( m_{2},s_{2} \right) \; \text{if} \; s\cdot
				s_2\cdot m_1 = s\cdot s_1\cdot m_2 \;\; \text{for some } s\in S.
				\label{eq:5.10}
\end{equation}
We denote by $m/s$ the equivalence class of $(m,s)$. \\
The functor $M\mapsto S^{-1}M$ commutes with colimits, finite-limits, tensor products and
Hom's. We obtain a functor 
\begin{equation}
				\begin{array}[h]{c}
								\begin{array}[h]{rl}
												\aset & \xhookrightarrow{\qquad} \bigo_{A}\text{-Set} \\\\
												M &\xmapsto{\qquad} M^{\natural} 
								\end{array}
\\\\
								M^{\natural}:= \left\{ \sigma\; :\; U\to \coprod_{\plusprimeideal\in U }
												S^{-1}_{\plusprimeideal}M,\; \sigma(\plusprimeideal)\in
												S^{-1}_{\plusprimeideal}M,\; \text{locally constant}
								\right\}
				\end{array}
				\label{eq:5.11}
\end{equation}
where $\sigma$ is {\underline{locally constant}} if for all $\plusprimeideal \in U$, there
is a neighborhood $\plusprimeideal\in D^{+}(s)\subseteq U$, and $m\in M$ such that for
$q\in D^{+}(s)$, $\sigma(q)\equiv m/s \in S_{q}^{-1}M$. \vspace{.1cm}\\
This functor prolongs naturally to functors 
\begin{equation}
				\begin{array}[H]{rl}
								\asetnab & \xhookrightarrow{\qquad} \bigo_{A}\text{-Set}^{\Bnabla} \\\\
								\Sigma(A) & \xhookrightarrow{\qquad} \Sigma(\bigo_{A}) \\\\
								S^{\cdot}_{A}\mymod & \xhookrightarrow{\qquad}
								S^{\cdot}_{\bigo_{A}}\mymod \\\\
								M & \xmapsto{\qquad} M^{\natural}
				\end{array}
				\label{eq:5.12}
\end{equation}
We have an identification of the stalk of $M^{\natural}$ at $\plusprimeideal\in
\text{spec}^{+}(A)$, 
\begin{equation}
				\begin{array}[H]{rl}
								M^{\natural}\Big|_{\plusprimeideal} =
								\lim\limits_{\xrightarrow[U\ni\plusprimeideal]{}} M^{\natural}(U)
								&\xrightarrow{\;\sim\;} S^{-1}_{\plusprimeideal}M := M_{\plusprimeideal}
								\\\\
								\left( \sigma,U \right)/\sim &\xmapsto{\quad\;} \sigma(\plusprimeideal)
				\end{array}
				\label{eq:5.13}
\end{equation}
Moreover, we have the identification of the global sections of $M^{\natural}$ over an affine basic
open subset $D^{+}(s)\subseteq\text{spec}^{+}(A)$, 
\begin{equation}
				\begin{array}[H]{rl}
								\Psi : M\left[ 1/s \right]:= \left\{ s^{\mathbbm{N}}\right\}^{-1}M 
								&\xrightarrow{\;\sim\;}M^{\natural}\left( D^{+}(s) \right) \\
								\frac{m}{s^n} &\xmapsto{\quad\;}\sigma(\plusprimeideal)\equiv
								\frac{m}{s^{n}}\in S^{-1}_{\plusprimeideal}M
				\end{array}
				 \label{eq: 5.14}
\end{equation}
\begin{proof}
				The map $\Psi$ which takes $m/s^{n}\in M\left[ 1/s \right]$ to the constant
				section $\sigma$ with $\sigma(\plusprimeideal)\equiv m/s^{n}$ for all
				$\plusprimeideal\in D^{+}(s)$ is well defined and is a map of $A[1/s]$-Sets. \\ 
				$\Psi$ {\underline{is injective}}: Assume $\Psi\left( m_{1}/s^{n_1}
				\right)=\Psi\left( m_{2}/s^{n_2} \right)$. Let $\plusideal$ denote the $+$ideal
				generated by $\left\{ a\in A^{+}_{[1]}, a\cdot s^{n_2}\cdot m_1=a\cdot
				s^{n_1}\cdot m_2 \right\}$. We have 
				\begin{equation}
								m_1/s^{n_1} = m_2/s^{n_2} \in A_{\plusprimeideal} \quad \text{for all
								$\plusprimeideal \in D^{+}(s)$}
								\label{eq:5.15}
				\end{equation}
				\begin{itemize}
								\item[$\Rightarrow$ ] $ s_{\plusprimeideal}\cdot s^{n_2}\cdot m_1 =
												s_{\plusprimeideal}\cdot s^{n_1}\cdot m_{2} $ with
												$s_\plusprimeideal\in A_{[1]}^{+}\setminus \plusprimeideal$, all
												$\plusprimeideal\in D^{+}(s)$
								\item[$\Rightarrow$ ] ${\plusideal}\not\subseteq\plusprimeideal$, all $\plusprimeideal\in D^{+}(s)$
								\item[$\Rightarrow$ ] $ V^{+}(\plusideal)\cap D^{+}(s)=\emptyset $
								\item[$\Rightarrow$ ] $ V^{+}(\plusideal)\subseteq V^{+}\left( s
												\right) $
								\item[$\Rightarrow$ ] $ s^{n}\in \plusideal $ for some $n>0$
								\item[$\Rightarrow$ ] $ s^{n}\cdot s^{n_2}_{2}\cdot m_{1}=  s^{n}\cdot s^{n_1}_{1}\cdot m_{2}$
								\item[$\Rightarrow$ ] $ m_1/s^{n_1} = m_2/s^{n_2} $ in $M\left[ 1/s
												\right]$
				\end{itemize}
				$\Psi$ {\underline{is surjective}}: Fix $\sigma\in M^{\natural}\left(
								D^{+}(s)
				\right)$. We get a covering by basic open subsets on which $\sigma$ is constant, 
				\begin{equation}
								D^{+}(s)= D^{+}(s_1)\cup \ldots \cup D^{+}(s_{N}), \quad
								\sigma(\plusprimeideal)\equiv m_i/s_i \;\; \text{for } \plusprimeideal\in
								D^{+}(s_i).
								\label{eq:5.16}
				\end{equation}
				On the subsets $D^{+}(s_i s_j)=D^{+}(s_i)\cap D^{+}(s_j)$ we have
				$m_i/s_i=m_j/s_j$, and by the injectivity of $\Psi: (s_is_j)^{n}\cdot s_j\cdot m_i
				= (s_i s_j)^{n}\cdot s_i \cdot m_j$ (with one $n$ which works for all $i,j\le N$).
				Replacing $s_{i}^{n}\cdot m_i$ by $m_i$, and $s_i^{n+1}$ by $s_i$, we may assume
				$\sigma$ is given by $m_i/s_i$ on $D^{+}(s_i)$, and moreover for all $i,j\le N: \;
				s_j\cdot m_i= s_{i}\cdot m_{j}$. Since $D^{+}(s)\subseteq \bigcup\limits_{i}
				D^{+}(s_i)$, we have for some $k>0$, $b,d\in A_{\ell}$, $j:\ell\to \left\{
								1,\ldots, N
				\right\}$,
				\begin{equation}
								s^{k}= \left( b\lideal s_{j(x)} \right)\sslash d
								\label{eq:5.17}
				\end{equation}
				Define $m\in M$ by 
			\begin{equation}
							m := \left< b,m_{j(x)}, d\right>
							\label{eq:5.18}
			\end{equation}
			We have for $i=1,\ldots, N$
			\begin{equation}
							\begin{array}[H]{lll}
											s_i\cdot m &= s_i\cdot \left< b,m_{j(x)},d\right> \\\\
											&= \left<b,s_{i} \cdot m_{j(x)},d\right> \\\\
											&= \left< b,s_{j(x)} \cdot  m_i,d \right> \\\\
											&= \left< b\lideal s_{j(x)},m_i,d \right> \\\\
											&= \left< (b\lideal s_{j(x)})\sslash d, m_{i},1 \right> \\\\
											&=  s^k\cdot m_i 
							\end{array}
							\label{eq:5.19}
			\end{equation}
			Thus $m_i/s_{i} = m/s^{k}$ in $M\left[ 1/s \right]$ and the section $\sigma$ is
			constant $\sigma=\Psi(m/s^k)$, and $\Psi$ is serjective. 
\end{proof}
Given a generalized scheme $X$, we have the full subcategories of ``quasi-coherent''
sheaves
\begin{equation}
				\begin{array}[h]{lrll}
								& q.c.\bigo_{X}\text{-Set} &\subseteq& \bigo_X\text{-Set} \\\\
								\text{resp., } & q.c.\bigo_{X}\text{-Set}^{\Bnabla} &\subseteq& \bigo_X\text{-Set}^{\Bnabla} \\\\
								& q.c.  \Sigma\left( \bigo_{X} \right) &\subseteq& \Sigma\left( \bigo_X  \right) \\\\
								& q.c. S^{\cdot}_{\bigo_X}\mymod &\subseteq& S^{\cdot}_{\bigo_X}\mymod
				\end{array}
				\label{eq:5.20}
\end{equation}
where an object $M$ is quasi-coherent if there is a covering of $X$ by open affines
$X=\cup_{i}\text{spec}^{+}(A_i)$, and there are $M_i\in A_i\text{-Set}$
(resp. $A_i\text{-Set}^{\Bnabla}$, $\Sigma(A_i)$, $S^{\cdot}_{A_{i}}\mymod$), with 
$M\big|_{\text{spec}^{+}(A_i)}\cong M_{i}^{\natural}$. Equivalently, for all open affine
$\text{spec}^{+}(A)\subseteq X$, we have 
\begin{equation}
M\big|_{{\text{spec}}^{+}(A)}\cong M\left( \text{spec}^{+}(A) \right)^{\natural}.
				\label{eq:5.21}
\end{equation}
For an affine scheme $X=\text{spec}^{+}(A)$ we get equivalences of categories
\begin{equation}
				\begin{array}[H]{l}
								\begin{array}[h]{l}
												\begin{tikzpicture}
																\node (Mnatural) at (0,4) {$M^{\natural}$};
																\node (M) at (0,2) {$M$};
																\draw [Bar->] (M) to (Mnatural);
																\node (A) at (4,4) {$q.c. \bigo_X\text{-Set}$};
																\node (B) at (4,2){$\aset$};
																\draw [->] (A) to [out=330,in=30,left]
																node{$\cong \hspace{0.6cm}$}
																(B);
																\draw [<-] (A) to [out=210,in=150,right]
																node{$\ $}
																(B);

																\node (A0) at (8,4) {$q.c. \bigo_{X}\text{-Set}^{\Bnabla}$};
																\node (B0) at (8,2){$\asetnab$};
																\draw [->] (A0) to [out=330,in=30,left]
																node{$\cong \hspace{0.6cm}$}
																(B0);
																\draw [<-] (A0) to [out=210,in=150,left]
																node{$\ $}
																(B0);
												\end{tikzpicture} 
								\end{array}
								\\\\
								\ \hspace{-0.8cm} \
								\begin{array}[h]{l}
												\begin{tikzpicture}
																\node (m) at (8,4) {$M$};
																\node (mx) at (8,2) {$M(X)$};
																\draw [Bar->] (m) to (mx);
																\node (A) at (0,4) {$q.c. \Sigma(\bigo_{X})$};
																\node (B) at (0,2){$\Sigma(A)$};
																\draw [->] (A) to [out=330,in=30,left]
																node{$\cong \hspace{0.6cm}$}
																(B);
																\draw [<-] (A) to [out=210,in=150,right]
																node{$\ $}
																(B);

																\node (A0) at (4,4) {$q.c.
																S^{\cdot}_{\bigo_{X}}\mymod$};
																\node (B0) at (4,2){$S^{\cdot}_{A}\mymod$};
																\draw [->] (A0) to [out=330,in=30,left]
																node{$\cong \hspace{0.6cm}$}
																(B0);
																\draw [<-] (A0) to [out=210,in=150,left]
																node{$\ $}
																(B0);
												\end{tikzpicture} 
								\end{array}
				\end{array}
				\label{eq:5.22}
\end{equation}
The subcategories of quasi-coherent sheaves (\ref{eq:5.20}) are closed under colimits and
finite-limits in the associated categories of sheaves, and are therefor co-complete and
finite-complete. They are also closed under the symmetric monoidal structure, so if
$M,N\in q.c. \bigo_X\text{-Set}$ (resp. $q. c.
				\bigo_{X}\text{-Set}^{\Bnabla}$; $q. c. \Sigma(\bigo_{X})$; $q. c.
S^{\cdot}_{\bigo_{X}}\mymod$) then so are $M\otimes_{\bigo_{X}}N$ (resp.
				$M\otimes_{\bigo_X}^{\Bnabla} N$; $M^{\cdot}\otimes_{\Sigma(\bigo_{X})}N^{\cdot}$;
$M^{\cdot}\otimes_{S^{\cdot}_{\bigo_X}}N^{\cdot}$) and $\myhom_{\bigo_{X}}(M,N)$
(resp.  $\myhom_{\bigo_{X}}^{\Bnabla}(M,N)$;
$\myhom_{\Sigma(\bigo_{X})}(M^{\cdot},N^{\cdot})$;
$\myhom_{S^{\cdot}_{\bigo_{X}}}\left( M^{\cdot},N^{\cdot} \right)$).
Moreover, for
any affine open subset $U=\text{spec}^{+}(A)\subseteq X$ we have
\begin{equation}
				\begin{array}[H]{rll}
								M\otimes_{\bigo_X}N\Big|_{U} &\cong  \left( M(U)\otimes_{A}N(U) \right)^{\natural} \\\\
								\myhom_{\bigo_{X}}(M,N)\Big|_{U} &\cong  \myhom_{A}\left( M(U),N(U) \right)^{\natural}
				\end{array}
				\label{eq:5.23}
\end{equation}
and similarly for $\otimes_{\bigo_{X}^{\Bnabla}}$, $\otimes_{\Sigma(\bigo_{X})}$,
$\otimes_{S^{\cdot}_{\bigo_{X}}}$, and $\myhom^{\Bnabla}_{\bigo_{X}}$,
$\myhom_{\Sigma(\bigo_{X})}$, $\myhom_{S^{\cdot}_{\bigo_{X}}}$.\vspace{.1cm}\\
For a mapping of generalized schemes $f:X\to Y$, the functor $f^{*}$ takes quasi-coherent
$\bigo_{Y}\text{-Sets}$ to quasi-coherent $\bigo_{X}\text{\normalfont-Set}$ and we get
strict-monoidal functors
\begin{equation}
				\begin{array}[H]{rll}
								f^{*}:q.c.\bigo_{Y}\text{-Set} &\xrightarrow{\qquad} q. c. \bigo_{X}\text{-Set} \\\\
								f^{*}:q.c.\bigo_{Y}\text{-Set}^{\Bnabla} &\xrightarrow{\qquad} q. c. \bigo_{X}\text{-Set}^{\Bnabla} \\\\
								f^{*}:q.c.\Sigma(\bigo_{Y}) &\xrightarrow{\qquad} q. c. \Sigma(\bigo_{X}) \\\\
								f^{*}:q.c.S^{\cdot}_{\bigo_{Y}}\mymod &\xrightarrow{\qquad}  q. c.
								S^{\cdot}_{\bigo_{X}}\mymod 
				\end{array}
				\label{eq:5.24}
\end{equation}
Moreover, for affine open subsets $U=\text{spec}^{+}(A)\subseteq X$,
$\mathcal{V}=\text{spec}^{+}(B)\subseteq Y$, with  $f(U)\subseteq \mathcal{V}$ and with
$\varphi=f^{\#}:B\to A$ the associated homomorphism, we have for quasi-coherent
sheaf $M$ on $Y$, 
\begin{equation}
				f^{*}M\Big|_{U} \cong \left( \varphi^{*} M(\mathcal{V}) \right)^{\natural}.
				\label{eq:5.25}
\end{equation}

								\begin{definition} \label{def:5.26}
												A mapping of generalized schemes $f:X\to Y$ is called a
												$c\text{\normalfont-map}$
												(``compact and quasi-separated'') if we can cover $Y$ by open affine
												subsets (or equivalently, if for all open affine)
												$U=\text{\normalfont spec}^{+}(B)\subseteq Y$, we have that $f^{-1}(U)$ is compact, so that
												\[
																f^{-1}(U)=\bigcup\limits_{i=1}^{N}\mathcal{V}_{i}, \qquad
																\mathcal{V}_{i}= \text{\normalfont spec}^{+}(A_i),
												\]
												is a \underline{finite} union of open affine subsets, and if moreover for
												$i,j\le N$, we have $\mathcal{V}_{i}\cap \mathcal{V}_{j}$ also compact, so that 
												\[
																\mathcal{V}_{i}\cap \mathcal{V}_j = \bigcup_{k=1}^{N_{i
																j}}\mathcal{W}_{i j k}, \qquad \mathcal{W}_{i j k} =
																\text{\normalfont spec}^{+}\left( A_{i j k} \right). 
												\]
\end{definition}
For a $c\text{\normalfont -map}$ $f:X\to Y$, the functor $f_{*}$ takes quasi-coherent sheaves
								to quasi-coherent sheaves, so we get adjunctions
								\begin{equation}
								\begin{array}[h]{l}
								\begin{tikzpicture}
												\node (A0) at (1,4) {$q.c.\bigo_{X}\text{-Set}$};
												\node (B0) at (1,2){$q.c.\bigo_{Y}\text{-Set}$};
												\draw [->] (A0) to [out=330,in=30,right]
												node{$f_{*}$}
												(B0);
												\draw [<-] (A0) to [out=210,in=150,left]
												node{$f^{*}$}
												(B0);

												\node (A) at (5,4) {$q.c.\bigo_{X}\text{-Set}^{\Bnabla}$};
												\node (B) at (5,2){$q.c.\bigo_{Y}\text{-Set}^{\Bnabla}$};
												\draw [->] (A) to [out=330,in=30,right]
												node{$f_{*}$}
												(B);
												\draw [<-] (A) to [out=210,in=150,left]
												node{$f^{*}$}
												(B);
								\end{tikzpicture}  \\\\
								\begin{tikzpicture}
												\node (A2) at (1,4) {$q.c.\Sigma (\bigo_{X})$};
												\node (B2) at (1,2){$q.c. \Sigma (\bigo_{Y})$};
												\draw [->] (A2) to [out=330,in=30,right]
												node{$f_{*}$}
												(B2);
												\draw [<-] (A2) to [out=210,in=150,left]
												node{$f^{*}$}
												(B2);

												\node (A3) at (5,4) {$q.c.S_{\bigo_{X}}^{\cdot}\mymod$};
												\node (B3) at (5,2){$q.c.S_{\bigo_{Y}}^{\cdot}\mymod$};
												\draw [->] (A3) to [out=330,in=30,right]
												node{$f_{*}$}
												(B3);
												\draw [<-] (A3) to [out=210,in=150,left]
												node{$f^{*}$}
												(B3);
								\end{tikzpicture}
				\end{array} 
				\label{eq:5.27}
								\end{equation}
								Indeed, for a quasi-coherent sheaf $M$ on $X$, we have 
								(using the above notations)  that 
								$f_{*}M|_U=f_{*}M|_{\text{spec}^{+}(B)}$ is the equalizer
								\begin{equation}
												\begin{array}[h]{l}
												        f_{*}M\big|_{U}\xrightarrow{\qquad}\prod\limits_{i=1}^{N}f_{*}\left(
												        M\big|_{\mathcal{V}_{i}} \right)\rightTwoArrows
												        \prod\limits_{i,j\le N}\prod\limits_{k=1}^{N_{i j}}f_{*}\left(
												        M\big|_{\mathcal{W}_{i j k}} \right) \\\\
																\text{and } 
																\begin{array}[t]{l}
																f_{*}\left( M\big|_{\mathcal{V}_i}
																\right) = \left( f_{i}^{\#}M(\mathcal{V}_{i})
																\right)^{\natural}, \\
																 f_{*}\left( M\big|_{\mathcal{W}_{i j k}}
																\right) = \left( f_{i j k}^{\#}M(\mathcal{W}_{i j k})
																\right)^{\natural}
																\end{array}
												\end{array}
												\label{eq:5.28}
								\end{equation}
								are quasi-coherent. \vspace{.1cm}\\
								For a generalized scheme $X$, let
								$\mathcal{D}(X)=\text{Aff}/X \ltimes S^\cdot\mymod$ denote the simplicial
								category with objects pairs $(\text{spec}^{+}(A),M^{\cdot})$, with
								$\text{spec}^{+}(A)\subseteq X$ open affine, and $M^{.}\in
								S^{\cdot}_{A}\mymod$. The maps are given by
								\begin{equation}
												\mathcal{D}(X)\left(
												(\text{\normalfont spec}^{+}(A_0),M_{0}^{\cdot}),
				(\text{\normalfont spec}^{+}(A_1),M_{1}^{\cdot})\right) =
				\R\text{Map}_{S^{\cdot}_{A_0}}\left(
				\mathbb{L}\phi^{*}M^{\cdot}_{1},M^{\cdot}_{0} \right)
				\label{eq:5.29}
								\end{equation}
								where
								$\phi:\text{spec}^{+}(A_0)\xhookrightarrow{\;\;}\text{spec}^{+}(A_1)$
								denotes the inclusion (and are empty if $\text{spec}^{+}\left( A_0
								\right)\not\subseteq\text{spec}^{+}(A_1)$). We have a forgetfull functor 
								\begin{equation}
												\rho: \mathcal{D}(X)\xrightarrow{\qquad} \text{Aff}/X. 
												\label{eq:5.30}
								\end{equation}
								The coherent nerve of $\mathcal{D}(X)$ is the simplicial set
								$N_{\BDelta}\left(\mathcal{D}(X)\right)\in \text{Set}^{\Bnabla}$ with
								$n$ simplices
								\begin{equation}
												N_{\BDelta}\left( \mathcal{D}(X) \right)_{n} =
												\text{Cat}_{\BDelta}\left(\mathcal{C}(\Delta^{n}),\mathcal{D}(X) \right)  
												\label{eq:5.31}
								\end{equation}
								Explicitly, its elements are given by the data of 
								\begin{equation}
												\begin{array}[H]{lll}
																\text{(\romannumeral 1)} &
																				(U_0\subseteq U_1 \subseteq
																				\dots \subseteq U_n)\in N_n\left( \text{Aff}/X \right), \quad
																				U_i=\text{spec}^{+}(A_i)\subseteq X \text{ open}, \\\\
																\text{(\romannumeral 2)} & 
																				M^{\cdot}_{i}\in S^{\cdot}_{A_{i}}\mymod \\\\
																\text{(\romannumeral 3)} &
																				\text{For $0\le i< j\le
																				n$}, \quad \varphi_{ij}\in
																				S^{\cdot}_{A_{i}}\mymod\left(
																				\mathbb{L}\phi_{i
																				j}^{*}M^{\cdot}_{j},M_{i}^{\cdot}
																				\right) 
												\end{array}
												\label{eq:5.32}
								\end{equation}
								where $\phi_{i j}: U_i\hookrightarrow U_j$ denotes the inclusion, and
								{\underline{higher coherent homotopies}} given by
								\begin{equation}
												\begin{array}[h]{l}
																\varphi_\tau\in S^{\cdot}_{A_i}\mymod\left(
																				\mathbb{L}\phi^{*}_{i j} M^{\cdot}_{j}\otimes
																\Delta(\ell)_{+},M_{i}^{\cdot} \right) \\\\
																\ \qquad\text{for all
																				$\tau=(\tau_0\subseteq \tau_1\subseteq \dots \subseteq \tau_\ell)\in
																N_{\ell}(\mathcal{P}_{i j})$}
												\end{array}
																\label{eq:5.53}
								\end{equation}
								(Where $P_{i j}$ is the set of subsets of $[i,j]$ containing the end
								points). \vspace{.1cm}\\
								These are required to satisfy:
								\begin{equation}
												\displaystyle d_{k}\varphi_{\tau}=
								\varphi_{(\tau_0\subseteq\dots
								\hat{\tau}_k\ldots \subseteq \tau_{\ell})}
								\label{eq:5.33}
								\end{equation} 
				and 
				\begin{equation}
								\begin{array}[h]{l}
												\text{for $\ell=0$, $\tau=(\tau_{0})=\left\{
												i,k_1,\dots,k_{q},j \right\}$}, \\\\
												 \varphi_{\tau_{0}} =  \varphi_{i k_1}\circ 
												 \mathbb{L}\phi_{i k_1}^{*}\left( \varphi_{k_1,k_2}\circ
																\mathbb{L}\phi_{k_1,k_2}^{*}\left( \ldots
																				\varphi_{k_{q-1},k_q}\circ \mathbb{L}
																				\phi_{k_{q-1},k_{q}}^{*}\left( \varphi_{k_q,j}
																\right)\dots
								\right) \right)
				\end{array}
								\label{eq:5.34}
				\end{equation}
				Forgetting all the data but the $n$ simplex $(U_0\subseteq \dots \subseteq U_n)$
				of $\text{Aff}/X$ gives an inner, Cartesian and co-Cartesian, fibration of
				quasi-categories 
				\begin{equation}
								N_{\BDelta}(\rho):
								N_{\BDelta}\left(\mathcal{D}(X)\right)\xtworightarrow{\qquad} N\left( \text{Aff}/X
								\right).
								\label{eq:5.35}
				\end{equation}
				The derived category of quasi-coherent symmetric spectra
				$\mathbb{D}(X)=\text{Ho}
				D(X)_{\infty}$ is the homotopy category of $D(X)_{\infty}$ the
				($\infty$-categorical) limit of $D(A)_{\infty}$,
				where $D(A)_{\infty}$ is
				the $\infty$-category of fibrant
				co-fibrant $S^{\cdot}_{A}\mymod$ associated to
				$\mathbb{D}(A)=\text{Ho}(S^{\cdot}_{A}\mymod)$, and the limit taken over all $A$ with
				$\text{spec}^{+}(A)\in \text{Aff}/X$. \vspace{.1cm}\\ 
				By Lurie's straightening/unstraightening
				functors, \cite{MR2522659}, Corollary 3.3.3.2, 
				this limit $D(X)_{\infty}$ is equivalent to the full sub-$\infty$-category on {\underline{Cartesian}}
				sections of the $\infty$-category of sections of (\ref{eq:5.35}) 
				\begin{equation}
								\begin{array}[h]{ll}
												\text{where a section } & m: N\left( \text{Aff}/X \right)\xrightarrow{\qquad} N_{\BDelta}\left(
																\mathcal{D}(X)
												\right) \\\\
												\ & N_{\BDelta }(\rho)\circ m =\text{id}
								\end{array}
								\label{eq:5.37}
				\end{equation}
				is {\underline{Cartesian}} iff for all $n$-simplex
				$(U_0\subseteq U_1\subseteq \dots \subseteq U_n)\in N_{n}\left( \text{Aff}/X \right)$,
				$U_i=\text{spec}^{+}(A_i)$, the data $m(U_{\cdot})$ in the notaions of
				(\ref{eq:5.32}) satisfies that all the arrows (\romannumeral 3)
				$\varphi_{i j}\in S^{\cdot}_{A_{i}}\mymod(\mathbb{L}\phi^{*}_{i j}
				M_{j}^{\cdot}, M_{i}^{\cdot})$, $0\le i < j\le n$, are weak-equivalences. The
				category $D(X)_{\infty}$ is symmetric monoidal stable quasi-category, so
				$\mathbb{D}(X)$ is symmetric monoidal triangulated category. 
												
\nocite{*}
\bibliography{algebra_over_generalized_rings}{}
\bibliographystyle{plain}
\end{document}